\documentstyle{amsart}

\newcommand{\nc}{\newcommand}
\nc{\C}{{\cal C}}
\nc{\W}{{\cal W}}
\renewcommand{\O}{{\cal O}}
\nc{\B}{{\cal B}}
\nc{\A}{{\frak A}}

\nc{\N}{{\mathbf N}}
\nc{\CC}{{\mathbf C}}
\nc{\Z}{{\mathbf Z}}
\nc{\R}{{\mathbf R}}
\nc{\si}{{\frac\infty 2}}
\nc{\St}{\operatorname{St}^{\frac\infty 2+\bullet}}
\nc{\gr}{\operatorname{gr}}
\renewcommand{\lim}{\operatorname{lim}}
\nc{\tr}{\operatorname{tr}}
\nc{\ch}{\operatorname{ch}}
\nc{\p}{{\frak p}}
\nc{\g}{{\frak g}}
\nc{\np}{{\frak n}^+}
\nc{\ap}{{\frak a}^+}
\nc{\am}{{\frak a}^-}
\nc{\nm}{{\frak n}^-}
\nc{\hatg}{\hat{\frak g}}
\renewcommand{\a}{{\frak a}}
\renewcommand{\b}{{\frak b}}
\nc{\gl}{{\frak g\frak l}}
\nc{\n}{{\frak n}}
\nc{\pl}{\partial}
\nc{\M}{{\cal M}}
\nc{\lth}{\ell t}
\nc{\cl}{{\cal C}\ell(\n\oplus\n^*)}
\nc{\h}{{\frak h}}
\nc{\overh}{\overline{{\frak h}}}
\nc{\overn}{\overline{{\frak n}}}
\nc{\overe}{\overline{e}}
\nc{\overa}{\overline{{\frak a}}}
\nc{\overb}{\overline{{\frak b}}}
\nc{\overr}{\overline{R}}
\nc{\hgt}{\operatorname{ht}}
\nc{\overg}{\overline{{\frak g}}}
\nc{\gothu}{{\frak u}}
\nc{\Ind}{\operatorname{Ind}}
\nc{\Coind}{\operatorname{Coind}}
\nc{\opp}{{\operatorname{opp}}}
\nc{\Ker}{\operatorname{Ker}}
\nc{\im}{\operatorname{Im}}
\nc{\Coker}{\operatorname{Coker}}
\nc{\inlim}{\underset{\map}{\operatorname{lim}}}

\nc{\Ext}{\operatorname{Ext}^{\bullet}}
\nc{\ext}{\operatorname{Ext}}
\nc{\ex}{\operatorname{exp}}
\nc{\Tor}{\operatorname{Tor}_{\bullet}}
\nc{\tor}{\operatorname{Tor}}
\nc{\Tors}{\operatorname{Tor}_{\frac \infty 2+\bullet}}
\nc{\Exts}{\operatorname{Ext}^{\frac \infty 2+\bullet}}
\nc{\Hom}{\operatorname{Hom}^{\bullet}}
\nc{\End}{\operatorname{End}}
\nc{\afr}{\{0\}\cup\underline{r}}

\nc{\ad}{\operatorname{ad}}
\renewcommand{\hom}{\operatorname{Hom}}

\renewcommand{\mod}{\operatorname{-mod}}
\nc{\Mod}{\operatorname{Mod}}

\nc{\Barb}{\operatorname{Bar}^{\bullet}}
\nc{\stan}{\underline{\Lambda(\n)}}
\nc{\cost}{\underline{\Lambda(\n^*)}}

\nc{\underc}{\underline{\CC}}
\nc{\underb}{\underline{U(\b)}}
\nc{\underr}{\underline{r}}
\nc{\undern}{\underline{n}}

\nc{\underA}{\underline{A}}
\nc{\underB}{\underline{B}}
\nc{\underk}{\underline{k}}
\nc{\map}{\longrightarrow}

\nc{\Q}{{\mathbf Q}}
\nc{\bs}{\bigskip\\}
\nc{\ms}{\smallskip\\}

\nc{\til}{\widetilde}

\nc{\Lemma}{{\bf Lemma:\ }}
\nc{\Theorem}{{\bf Theorem:}\ }
\nc{\Cor}{{\bf Corollary:}\ }
\nc{\Def}{{\bf Definition:}\ }
\nc{\Prop}{{\bf Proposition:}\ }
\nc{\Con}{{\bf Conjecture:}\ }
\nc{\Rem}{{\bf Remark:}\ }
\nc{\dok}{{\bf Proof.}\ }

\nc{\bul}{^{\bullet}}
\nc{\Cone}{\operatorname{Cone}}
\nc{\supp}{\operatorname{supp}}
\nc{\ten}{\otimes}
\nc{\stand}{C^{\frac \infty 2+\bullet }}

\nc{\ssn}{\subsection{}}
\nc{\sssn}{\subsubsection{}}

\nc{\sqbinom}{\fracwithdelims[][0pt]}

\author{S.~M.~Arkhipov}
\title{A new construction of the semi-infinite BGG resolution}
\address{Independent University of Moscow, Pervomajskaya st. 16-18,
Moscow 105037, Russia}
\email{hippie@@ium.ips.ras.ru}
\begin{document}
\maketitle
\section{Introduction}
In their beautiful work [BGG] I.~N.~Bernshtein, I.~M.~Gelfand and
S.~I.~Gelfand introduced so called Bernshtein-Gelfand-Gelfand resolutions of
finite dimensional simple modules over semisimple Lie algebras consisting of
Verma modules. Some years after that similar resolutions of integrable simple
modules over general symmetrizable Kac-Moody algebras were constructed in
[RW]. In both cases the resolutions, though obtained by purely algebraic
methods, encoded geometric  structure of the Flag manifolds of the
corresponding Lie groups (the decomposition into the union of the Schubert
cells).

In [FF] B.~L.~Feigin and E.~Frenkel introduced some analogues of BGG
resolutions of integrable simple modules over affine Kac-Moody
algebras that consisted of Wakimoto modules. Mathematical folklore
says that Wakimoto modules and Feigin-Frenkel resolutions represent
some conjectural geometry of semi-infinite flag manifolds (see [FF]).
This geometry should include in particular a semi-infinite analogue of
the Schubert decomposition and a semi-infinite localization of a certain
category of modules over the affine Lie algebra.

\ssn The present paper is devoted to an attempt of better understanding of
the Feigin-Frenkel construction avoiding bosonization techniques and physical
vocabulary used in [FF]. We construct a sequence of analogues of the affine
BGG resolutions enumerated by elements of the affine Weyl group. Each of
these complexes is quasiisomorphic to the integrable simple module over the
affine Lie algebra. Moreover the complexes form an inductive system
enumerated by the dominant chamber elements of the corresponding weight
lattice and the limit of the system happens to be the semi-infinite BGG
resolution of Feigin and Frenkel.

I hope that the techniques of the present paper will allow to
construct semi-infinite BGG resolutions for general symmetrizable
Kac-Moody algebras. It it possible also to perform an analogue of our
construction for affine quantum groups for general $q$.  This will be
explained elsewhere.

\ssn
Let us describe the structure of the paper. In the second section we
develop the techniques of semiregular modules over Lie algebras.
Namely let $\g$ be a graded Lie algebra with a negatively graded
finite dimensional subalgebra $\n$.  The $\g$-module
$S_{\n}:=\Ind_{\n}^{\g}\Coind_{\CC}^{\n}\CC$ is called the semiregular
module over $\g$.  In the second and the third sections we present two
different and independent proofs of the fact that the endomorphism
algebra of $S_{\n}$ contains the universal enveloping algebra of $\g$.
Thus $S_{\n}$ becomes a $\g$-bimodule. This allows us to construct a
kernel functor $*\ten_{U(\g)}S_{\n}:\ \g\mod\map\g\mod$.

In the fourth section we recall several facts from combinatorics of
the affine root systems including definitions of the semi-infinite
Bruhat order and the semi-infinite length function.

In the fifth section we apply the techniques from the second and the
third sections in the case of affine Lie algebras. The twisting
functors $*\ten_{U(\g)}S_{\n_w}$ for some finite dimensional nilpotent
subalgebras $\n_w$ enumerated by  elements of the affine Weyl group
are the affine analogues of the twisting functors defined in [FF] in
the case of finite dimensional semisimple Lie algebras.  We call the
images of Verma modules under these functors the twisted Verma
modules.  Taking the images of the affine BGG resolutions under the
twisting functors we obtain a family of complexes enumerated by
elements of the affine Weyl group.  We call these complexes the
twisted BGG resolutions.

In the sixth section we show that twisted BGG resolutions form an
inductive system and when the lengths of elements of the Weyl group
tend to infinity, the corresponding twisted BGG resolutions ``tend to
the semi-infinite BGG resolution of Feigin and Frenkel''.

\section{Endomorphism algebra of the semiregular module}

\subsection{Calculations in Clifford algebra.}
Let $\n=\underset {m\in\\Z_{<0}}{\bigoplus}\n_m$ be a finite dimensional
graded vector space, $\dim\n=n$. Then there is a canonical bilinear form
on the vector space $\n\oplus\n^*$:
$$
\langle
(n_1,f_1),(n_2,f_2)\rangle=f_2(n_1)+f_1(n_2),\ n_1,n_2\in\n, \
f_1,f_2\in\n^*.
$$
\sssn
Recall that the {\em Clifford algebra} $\cl
$ is a $\CC$-algrbra with the space of
generators equal to $\n\oplus\n^*$ and the following relations:
$$
\omega_1\omega_2+\omega_2\omega_1=\langle\omega_1,\omega_2\rangle,\
\omega_1,\omega_2\in\n\oplus\n^*.
$$
\Rem
The subalgebra in $\cl$ generated by the space $\n$
(resp. the space $\n^*$) is isomorphic to the exterior algebra
$\Lambda(\n)$ (resp. the exterior algebra $\Lambda(\n^*)$).

\sssn
\Lemma
$\cl\cong \Lambda(\n)\ten\Lambda(\n^*)$ as a vector space.\qed

In particular $\dim\cl=2^{2(\dim\n)}$. From now on the set
$\{1,\ldots,n\}$ is denoted by $\underline{n}$.  Choose a homogeneous base
$\{e_i|i\in\underline{n}\}$ of the vector space $\n$ and its dual base
$\{e_i^*|i\in\underline{n}\}\subset\n^*$. Fix a subset
$I\subset\{1,\ldots,\dim\n\}$ and denote its complementary set by
$\overline{I}$. Denote the element of the Clifford algebra $\underset{i\in
I}{\prod} e_i$ (resp. the element $\underset{i\in I}{\prod} e_i^*$ ) by
$e_I$ (resp. by $e_I^*$).
The order of factors in the monomials is the one in the set
$\underline{n}$.

\sssn
Consider the left $\cl$-modules
$$
\stan:=\Ind_{\Lambda(\n^*)}^{\cl}\underc\text{ and }\cost:=
\Ind_{\Lambda(\n)}^{\cl}\underc.
$$
Here $\underc$ denotes the trivial module.
Then the images of the elements $\{e_I|I\subset\undern\}$
(resp. of the elements $\{e_I^*|I\subset\underline{n}\}$) form a base in the
space $\stan$ (resp. $\cost$). We are going to calculate the action of the
Clifford algebra in these bases explicitly. Denote the matrix element
$e_I\mapsto e_J$ (resp. $e_I^*\mapsto e_J^*$) by $a_I^J$ (resp. by
$b_I^J$).

Consider the elements in the Clifford algebra of the form
$
e_Je^*_{\underline{n}}e_I,\ I,J\subset\underline{n}.
$

\sssn
\Lemma                     \label{matr}

\qquad(i) The element $e_Je^*_{\underline{n}}e_I$ acts as $a_{\overline{I}}^
{J}:\ \stan\map\stan$;

\qquad(ii) The element $e_Je^*_{\underline{n}}e_I$ acts as
$b_I^{\overline{J}}:\ \cost\map\cost$.
\qed

\sssn
\Cor
In particular $\cl$ is isomorphic to the $(2^n\times 2^n)$ matrix
algebra.\qed

\ssn
The antiautomorphism of the tensor algebra $T(\n\oplus\n^*)$
$$
\sigma:\ a_1\ten\ldots\ten a_k\mapsto\a_k\ten\ldots\ten a_1,\
a_1,\ldots, a_k\in\n\oplus\n^*.
$$
preserves the relation defining the Clifford algebra, hence
the restriction of $\sigma$,
$$
\sigma:\ \cl\map\cl^{\opp},
$$
is correctly defined. The antiautomorphism $\sigma$ preserves the
subalgebras $\Lambda(\n)$ and $\Lambda(\n^*)$.
For any $\omega\in \Lambda^k(\n)$ or $\omega\in\Lambda^k(\n^*)$ we have
$\sigma(\omega)=(-1)^{\left[\frac k2\right]}\omega$, where $[\cdot]$
denotes the integral part of the number.

\sssn
\Rem
There is a natural isomorphism of algebras
$$
\beta:\ \End_{\CC}(\stan)\til{\map}\End_{\CC}(\cost)^{\opp},\
a_I^J\mapsto b_J^I,\ I,J\subset\underline{n}.
$$
Denote the isomorphism $\cl\til{\map}\End_{\CC}(\stan)$ by $\alpha$,
the isomorphism $\End_{\CC}(\cost)\til{\map}\cl$ is denoted by $\gamma$.

\sssn
\Lemma  \label{ident}
$\gamma\circ\beta\circ\alpha=\sigma .$

\dok
Follows immediately from~\ref{matr}. \qed

\subsection{Endomorphism algebra of the semiregular module}\label{n}
Suppose we have a graded Lie algebra $\g=\underset{m\in\Z}{\bigoplus}\g_m$
and a negatively graded finite dimensional Lie subalgebra $\n\subset\g,\
\dim\n=n$.
Note that $\n$ is automatically nilpotent. Suppose also that $\n$
acts locally ad-nilpotently on $\g$ (and thus on $U(\g)$ as well).
Fix a graded vector subspace
$\b\subset\g$ such that $\g=\b\oplus\n$ as a vector space. Denote
the $\g$-module $\Ind_{U(\n)}^{U(\g)}\underc$ by $\underb$.
Choose a homogeneous base
$\{e_i\}_{i\in\underline{n}}$ of $\n$. Thus the bracket in $\n$
is defined by the structure constants:
$
[e_i,e_j]=\sum_{k=1}^nc_{ij}^ke_k.
$

\sssn
The category of graded left $\g$-modules $M=\underset{m\in\Z}{\bigoplus}
M_m$ with morphisms being $\g$-module homomorphisms that preserve gradings
is denoted by $\g\mod$.  The corresponding category of $\n$-modules is denoted
by $\n$-mod. Let $\langle\cdot\rangle$ denote the grading shift functor:
$$
\langle\cdot\rangle:\ \g\mod\map\g\mod,\ M\in\g\mod,\text{ then }
M\langle k\rangle_{m}:=M_{k+m}.
$$
Denote the space $\underset{k\in\Z}{\bigoplus}
\hom_{\g\mod}(M_1,M_2\langle k\rangle)$
by $\hom_{U(\g)}(M_1,M_2)$.
A similar notation is used in the category of $\n$-modules.

\sssn
The grading on $\n$ induces a natural grading on its universal enveloping
algebra $U(\n)$.  Consider a {\em left} coregular $\n$-module
$U(\n)^*\in\n\mod$:
$$
U(\n)^*:=\underset{m\in\N}{\bigoplus}\hom_{\CC}(U(\n)_{-m},\CC).
$$
The left
action of $\n$ is defined as follows:
$$
f:\ U(\n)\map\CC,\ g\in\n,\text{
then } (n\cdot f)(u):=f(un),\ u\in U(\n).
$$
\Def
The left $\g$-module
$S_{\n}:=\Ind_{U(\n)}^{U(\g)}U(\n)^*$ is called the left semiregular
representation of $\g$ with respect to the subalgebra $\n$.

Our main task in this section is to find the endomorphism algebra
$\End_{U(\g)}(S_{\n})$. Let us calculate it first as a vector space.

\sssn
\Lemma
$\ext_{\g\mod}^{>0}(S_{\n},S_{\n})=0,\ \End_{U(\g)}(S_{\n})\cong\hom_{\CC}
(U(\n)^*,\underb)$.

\dok
By Schapiro lemma $\Ext_{U(\g)}(S_{\n},S_{\n})=\Ext_{U(\n)}(U(\n)^*,S_{\n})$.
But the  $U(\n)$-module $S_{\n}$ is cofree with the space of cogenerators
equal to $\underb$. Again use Schapiro lemma.  \qed

\sssn
\Rem Thus we see that up to a completion the space $\End_{U(\g)}(S_{\n})$
is equal to $U(\n)\ten\underb\cong U(\g)$.

The  considerations below remain true also for subalgebras $\n\subset\g$
such that
the adjoint action of $\n$ on $\g$ is not locally nilpotent. In this case
the algebra that is being calculated equals
to the endomorphism algebra of some completion of $S_{\n}$.

\ssn
Consider the $U(\g)$-free resolution  $P\bul(S_{\n})$ of $S_{\n}$ as follows:
\begin{gather*}
P^{-k}(S_{\n}):=U(\g)\ten\Lambda^{k}(\n)\ten U(\n)^*,\ d:\
P^{-k}(S_{\n})\map P^{-k+1}(S_{\n}),\\ u\in U(\g), f\in U(\n)^*,\
n_1,\ldots,n_k\in\n,\
d(u\ten n_1\wedge\ldots\wedge n_k\ten f):=\\
\sum_{j=1}^k(-1)^j(un_j\ten
n_1\wedge\ldots\wedge n_{j-1}\wedge n_{j+1}\wedge\ldots\wedge n_k\ten f+
u\ten
n_1\wedge\ldots\wedge n_{j-1}\wedge n_{j+1}\wedge\ldots\wedge
n_k\ten (n_j\cdot f))\\
+
\sum_{i<j}(-1)^{i+j}u\ten[n_i,n_j]\wedge n_1\wedge\ldots\wedge
n_{i-1}\wedge n_{i+1}\wedge\ldots\wedge n_{j-1}\wedge
n_{j+1}\wedge\ldots\wedge n_k\ten f.  \end{gather*}

\sssn Consider a graded subalgebra $A\bul\subset
\End_{U(\g)}(P\bul(S_{\n}))$:
$$
A\bul:=U(\g)\ten\End_{\CC}(\Lambda(\n)\ten U(\n)^*)
\til{\map}U(\g)\ten\cl\ten\End_{\CC}(U(\n)^*).
$$
\sssn
\Lemma
$A\bul$ is a DG-subalgebra in
$\End_{U(\g)}(P\bul(S_{\n}))$.\qed

Denote the elements of the chosen base of $\n$ (resp. of $\n^*$)
considered as elements of $\cl$
by $\{\overline{e}_i\}_{i\in\underline{n}}$, (resp.
by $\{\overline{e}^*_i\}_{i\in\underline{n}}$).
Consider the element in $A\bul$ as follows:
$$
D_{A\bul}:=\sum_{j=1}^k(e_i\ten \overline{e}^*_i\ten 1+
1\ten \overline{e}_i^*\ten l_{e_i})+
\sum_{i<j,k}c_{i,j}^k\overline{e}_i^*\overline{e}_j^*\overline{e}_k,
$$
where $l_{e_i}(f)$ denotes the element $(e_i\cdot f)\in U(\n)^*$.

\sssn
\Lemma The differential in the DG-algebra $A\bul$ is given by the
supercommutator
with the element $D_{A\bul}$: for every $a\in A\bul$ we have
$
d(a)=\{D_{A\bul},a\}.
$
\qed

\sssn
\Lemma
$H^{\ne 0}(A\bul)=0$. \qed

\ssn
We construct another DG-algebra of the same size. Consider the left regular
$U(\g)$-module $U(\g)$ and its $U(\g)$-free resolution $Q\bul(U(\g))$:
\begin{gather*}
Q^{-k}(U(\g)):=U(\g)\ten U(\n)\ten\Lambda^{k}(\n),\ d:\
Q^{-k}(U(\g))\map Q^{-k+1}(U(\g)),\\ u\in U(\g), m\in U(\n),\
n_1,\ldots,n_k\in\n,\
d(u\ten m\ten n_1\wedge\ldots\wedge n_k):=\\
\sum_{j=1}^k(-1)^j(un_j\ten m\ten
n_1\wedge\ldots\wedge n_{j-1}\wedge n_{j+1}\wedge\ldots\wedge n_k+
u\ten mn_j\ten
n_1\wedge\ldots\wedge n_{j-1}\wedge n_{j+1}\wedge\ldots\wedge
n_k)\\
+
\sum_{i<j}(-1)^{i+j}u\ten m\ten[n_i,n_j]
\wedge n_1\wedge\ldots\wedge n_{i-1}\wedge
n_{i+1}\wedge\ldots\wedge n_{j-1}\wedge n_{j+1}\wedge\ldots\wedge n_k
\end{gather*}
(this seems to be unnessesary since $U(\g)$ is free over itself).

\sssn
Consider a graded subalgebra $B\bul\subset\End_{U(\g)}(Q\bul(U(\g)))$:
$$
B\bul:=U(\g)\ten\End_{\CC}(\Lambda(\n)\ten U(\n))
\til{\map}U(\g)\ten\cl\ten\End_{\CC}(U(\n)).
$$
and an element $D_{B\bul}$  in $B\bul$ as follows:
$$
D_{B\bul}:=
\sum_{j=1}^k(e_i\ten \overline{e}^*_i\ten 1+
1\ten \overline{e}_i^*\ten r_{e_i})+
\sum_{i<j,k}c_{i,j}^k\overline{e}_i^*\overline{e}_j^*\overline{e}_k,
$$
where $r_{e_i}(m)$ denotes the element $(m\cdot e_i)\in U(\n)$.
\sssn
\Lemma
The differential in the DG-algebra $B\bul$ is given by the
supercommutator with the element $D_{B\bul}$: for every $b\in B\bul$ we have
$ d(b)=\{D_{B\bul},b\}.
$ \qed

\sssn
\Lemma
$
H^{\ne0}(B\bul)=0.
$
\qed

We construct an inclusion $\iota$ of the algebra $U(\g)$ into $B^0$ such
that $\iota:\ U(\g)\hookrightarrow B\bul$ is a quasiisomorphism of DG-algebras.

\sssn
\Lemma
For any right $\n$-module $X$ consider the right $\n$-module $X_1:=X\ten
U(\n)$ with $\n$-action provided by the Hopf algebra structure on $U(\n)$ and
the right $\n$-module $X_2:=X\ten U(\n)$ with $\n$-action trivial on the
first factor.  Then $X_1\til{\map}X_2$.

\dok
Consider the comultiplication map in the algebra $U(\n)$:
$$
\Delta:\ U(\n)\map U(\n)\ten U(\n),\ \Delta(u)=\sum_i\Delta^{(1)}_i(u)\ten
\Delta^{(2)}_i(u)+ 1\ten u,\ u\in U(\n).
$$
We introduce a map $\phi:\ X_2\map X_1$ as follows:
$$
\phi:\ x\ten u \mapsto x\ten u +\sum_i(x\cdot
\Delta^{(1)}_i(u))\ten\Delta^{(2)}_i(u),\ x\in X,\ u\in U(\n).
$$
Here $\cdot$ denotes the right $\n$-action on $X$. Evidently
$\phi$ is an isomorphism of vector spaces. One checks directly that
$\phi$ commutes with the defined $\n$-actions on $X_1$ and $X_2$.
\qed

Let $\widetilde{U(\g)\ten U(\n)}$ be the $U(\g)$-$U(\n)$ bimodule
with left $U(\g)$-action via the first factor and  right $U(\n)$-action
via the second one. Consider also the $U(\g)$-$U(\n)$ bimodule
$\overline{U(\g)\ten U(\n)}$ with left $U(\g)$-action via the first factor
and right $U(\n)$-action using the comultiplication on $U(\n)$.

\sssn
\Cor
$\phi:\ \til{U(\g)\ten U(\n)}\til{\map} \overline{U(\g)\ten U(\n)}$.
\qed

\Rem
Note that $Q\bul(U(\g))$ is the standard complex for the computation of the
Lie algebra homology of $\n$ with coefficients in the right module
$\overline{U(\g)\ten U(\n)}$.

Note also that the right $U(\g)$-action on $\til{U(\g)\ten U(\n)}$ commutes
with both the left $U(\g)$-module and the right $U(\n)$-module structures.
Thus we obtain the following statement.

\sssn
\Lemma
The map $g\mapsto \phi\circ (r_g\ten 1\ten 1)\circ\phi^{-1}$,
where $r_g$ denotes the right multiplication by $g$,
defines a DG-algebra homomorphism $U(\g)\map\End_{U(\g)}(U(\g)\ten U(\n)\ten
\Lambda\bul(\n))=B\bul$ that becomes an isomorphism on cohomologies.

\dok
Follows immediately from the previous Remark.
\qed

\sssn
\Rem
Clearly the map from the previous Lemma is a composition
$$
U(\g)\overset{\iota}{\map} B\bul\map\End_{U(\g)}(Q\bul(U(\g))).
$$
Our main result in this section is the following statement. Consider the
following isomorphisms
of graded algebras (not nessesarily preserving differentials):
\begin{gather*}
\theta:\
B^{\bullet\opp}=\left(U(\g)\ten\cl\ten\End_{\CC}(U(\n))\right)^{\opp}
\til{\map} U(\g)\ten\cl\ten\left(\End_{\CC}(U(\n))\right)^{\opp}\\
\theta=\sigma_{\g}\ten\sigma_{\n},\ \sigma_{\g}(g)=-g,\ g\in\g,\
\sigma(\overline{e}_i) =\overline{e}_i,\
\sigma(\overline{e}_i^*)=\overline{e}_i^*\,\ i\in\underline{n}\\ \eta:\
U(\g)\ten\cl\ten\left(\End_{\CC}(U(\n))\right)^{\opp}
\til{\map}U(\g)\ten\cl\ten\End_{\CC}(U(\n)^*)=A\bul.
\end{gather*}
\sssn
\Theorem   \label{main}
The map $\eta\circ\theta:\ B^{\bullet\opp}\map A\bul$ is homomorphism of
DG-algebras.

\Rem
In particular there is a canonical inclusion of algebras
$U(\g)\map
\End_{U(\g)}(S_{\n})$ that is an isomorphism up to a certain completion
of $\im U(\g)$.

\subsection{Proof of the Theorem~\ref{main}.}
We begin the proof with presenting several Lemmas.

\sssn                     \label{rem}
\Rem
Let $\n$ be an arbitrary  Lie algebra with a base
$\{e_i\}_{i\in\overline{n}}$
and structure constants $c_{i,j}^k,\ i,j,k\in\underline{n}$. Then for
a fixed $k\in\underline{n}$ we have
$
\sum_{i,j\in\overline{n}}c_{i,j}^k=0.
$

Let $D_{A\bul}^{(1)}$ and $D_{B\bul}^{(1)}$ be
 the  summands of $D_{A\bul}$ and $D_{B\bul}$ respectively that belong
 to $U(\g)\ten\cl$ as follows:
$$
D_{A\bul}^{(1)}:=
\sum_{j=1}^ke_i\ten \overline{e}^*_i\ten 1+
\sum_{i<j,k}c_{i,j}^k\overline{e}_i^*\overline{e}_j^*\overline{e}_k,
$$
and $D_{B\bul}^{(1)}\in B\bul$ given by the same formula.

\sssn
\Lemma
$$
\sigma(D_{B\bul}^{(1)})=
-\sum_{j=1}^ke_i\ten \overline{e}^*_i\ten 1-
\sum_{i<j,k}c_{i,j}^k\overline{e}_i^*\overline{e}_j^*\overline{e}_k,
$$
in particular $d^{(1)}:\ U(\g)\ten\cl\map U(\g)\ten\cl$, $d^{(1)}(\cdot):=
\{D_{B\bul}^{(1)},\cdot\}$, commutes with $\sigma$.

\dok
The only thing  to be checked is that
$$
\sigma_{\n}\left(\sum_{i<j,k}c_{i,j}^k\overline{e}_i^*\overline{e}_j^*
\overline{e}_k\right)=
-\sum_{i<j,k}c_{i,j}^k\overline{e}_i^*\overline{e}_j^*\overline{e}_k.
$$
The calculation looks as follows:
\begin{gather*}
\sigma_{\n}\left(\sum_{i<j,k}c_{i,j}^k\overline{e}_i^*
\overline{e}_j^*\overline{e}_k\right)=
\sum_{i<j,k}c_{i,j}^k\overline{e}_k
\overline{e}_j^*\overline{e}_i^*=
-\sum_{i<j,k}c_{i,j}^k\overline{e}_j^*
\overline{e}_k\overline{e}_i^*+
\sum_{i<j}c_{i,j}^j\overline{e}_i^*=\\
-\sum_{i<j,k}c_{i,j}^k\overline{e}_i^*\overline{e}_j^*\overline{e}_k+
\sum_{i<j}c_{i,j}^j\overline{e}_i^*-
\sum_{i<j}c_{i,j}^i\overline{e}_j^*=\\
-\sum_{i<j,k}c_{i,j}^k\overline{e}_i^*\overline{e}_j^*\overline{e}_k+
\sum_{i<j}c_{i,j}^j\overline{e}_i^*+
\sum_{i>j}c_{i,j}^j\overline{e}_i^*.
\end{gather*} The statement follows from Remark~\ref{rem}.
\qed

Let $D_{A\bul}^{(2)}:=D_{A\bul}-D_{A\bul}^{(1)}$
and $D_{B\bul}^{(2)}:=D_{B\bul}-D_{B\bul}^{(1)}$.
It remains to prove that $\eta\circ\theta(D_{B\bul}^{(2)})=-D_{A\bul}^{(2)}$.

For a left $\n$-module $X$ consider the $\cl\ten\End_{\CC}(X)$-modules
$\stan\ten X$ and $\cost\ten X$. The corresponding identifications
$$
\cl\ten\End_{\CC}(X)\cong\End_{\CC}(\stan\ten X)\text{ and }
\cl\ten\End_{\CC}(X)\cong\End_{\CC}(\cost\ten X)
$$
are denoted by $\varphi$ and $\psi$. Consider the map
$
d:\stan\ten X\map\stan\ten X,$
$$
a_1\wedge\ldots\wedge a_k\ten x\mapsto\sum_{j=1}^k(-1)^j
a_1\wedge\ldots\wedge a_{j-1}\wedge a_{j+1}\wedge\ldots
\wedge a_k\ten a_j\cdot x,
$$
(resp. the map
$
\overline{d}^*:\cost\ten X\map\cost\ten X,\
a_1^*\wedge\ldots\wedge a_k^*\ten x\mapsto\sum_{s\in\underline{n}}
\overline{e}_s^*
\wedge a_1^*\wedge\ldots\wedge a_k^*\ten e_s\cdot x).
$
\sssn
\Lemma
$\varphi^{-1}(d)=\psi^{-1}(\overline{d}^*)$.

\dok
Follows from~\ref{matr}
\qed

Now by~\ref{ident} the map $\eta\circ\theta$ can be viewed as follows:
\begin{gather*}
B^{\bullet\opp}=\left(U(\g)\ten\End_{\CC}(\stan\ten U(\n))\right)^{\opp}
\til{\map}U(\g)\ten\left(\End_{\CC}(\stan\ten U(\n))\right)^{\opp}\\
\til{\map}U(\g)\ten\End_{\CC}(\cost\ten U(\n)^*)
\til{\map}U(\g)\ten\End_{\CC}(\stan\ten U(\n)^*)=A\bul.
\end{gather*}
Here the map
$U(\g)\ten\End_{\CC}(\cost\ten U(\n)^*)
\til{\map}U(\g)\ten\End_{\CC}(\stan\ten U(\n)^*)$
is given by $\varphi\circ\psi^{-1}$ for $X=U(\n)^*$.
Note that the map $d$ in the previous Lemma is given by the supercommutator
with $D_{A\bul}^{(2)}$, and the map $\overline{d}^*$ is given by
the supercommutator with the image of $D_{B\bul}^{(2)}$. So the theorem is
proved.  \qed

\Rem
In fact we did not use the Lie algebra structure on $\g$ in the proof of
Theorem ~\ref{main}. So the statement of the Theorem holds in the following
situation.  Let $C$ be a graded associative algebra containing a finite
dimensional negatively graded nilpotent Lie subalgebra $\n$ that acts locally
ad-nilpotently on $C$.  Suppose that $C$ contains $U(\n)$ as a graded
subalgebra and is free both as a left and  a right $U(\n)$-module.  Suppose
also  there exists an antipode map $\alpha: C\map C^{\opp}$ that preserves
$\n\subset C$.  Consider the left semiregular $C$-module
$S_{\n}:=C\ten_{U(\n)}U(\n)^*$. Then there exists an inclusion of associative
algebras $$ C\map\End_C(S_{\n}).  $$

\section{Endomorphisms of the semiregular module (continued)}

In this section we present another proof of Theorem~\ref{main}. Like in the
previous section let $\g=\underset{m\in\Z}{\bigoplus}\g_m$ be a graded Lie
algebra with a finite dimensional negatively graded Lie subalgebra
$\n\subset\g$ that acts locally ad-nilpotently both on $\g$ and on $U(\g)$.

\ssn
Suppose $\n$ is abelian, and $\dim\n=n$.  Choose a homogenious base
$\{e_i|i\in\undern\}$ in $\n$. Let $\{x_i|i\in\undern\}$ be the dual base in
$\n^*$.  Then $U(\n)=\CC[e_i|i\in\undern]$ and
$U(\n)^*=\CC[x_i,i\in\undern]$.  The Lie algebra $\n$ acts on $U(\n)^*$ by
derivations: $e_i$ acts as $\pl/ \pl x_i$.

\sssn                 \label{exp}
Consider the left $\g\oplus\n$-module $U(\g)\ten U(\n)^*\cong
U(\g)\ten\CC[x_i| i\in\undern]$.  For every $u\in U(\g)$ consider an
endomorphism $\sigma(u)$ of $U(\g)\ten U(\n)^*$ as follows:
$$
\sigma(u)(v\ten p):=v\ex\left(\sum_{i=1}^n \ad_{e_i}\ten x_i\right)(u)p,\
v\in U(\g),\ p\in\CC[x_i|i\in\undern].
$$
\sssn
\Lemma
\label{act}
For every $u\in U(\g)$ the endomorphism $\sigma(u)$ commutes both
with the diagonal action of $U(\n)$ and with the left regular action of
$U(\g)$ on $U(\g)\ten U(\n)^*$.

\dok
A direct calculation.\qed

\Cor $\sigma$ defines an inclusion of algebras
$U(\g)\hookrightarrow\End_{\g}(S_{\n})$.  \qed

\ssn                 \label{filtration}
Let $\g$ and  $\n$ be like in the beginning of this section.
Then there exists a filtration $F$ on $\n$:
$$
\n=F^0\n\supset F^1\n\supset\ldots\supset F^{\operatorname{top}}\n=0
$$
such that each $F^i\n$ is an ideal in $F^{i-1}\n$ and there exist abelian
subalgebras $\n^i\subset F^i\n$ such that each $F^i\n=\n^i\oplus F^{i+1}\n$
as a vector space.

\sssn
\Theorem                                               \label{main2}
There exists  an inclusion of algebras
$U(\g)\hookrightarrow\End_{\g}(S_{\n})$ such that the image of $U(\g)$ is
a dense subalgebra in $\End_{\g}(S_{\n})$.

\sssn
\Rem
The rest of this subsection is devoted to the proof of Theorem~\ref{main2}.

Suppose $\g\supset\n=\np\oplus\nm$ as a graded vector space, where $\nm$ is
an ideal in $\n$ and $\np$ is a subalgebra in $\n$.  Consider the left
$\n$-modules
$$
S_{\nm}^{\n}:=U(\n)\ten_{U(\nm)}U(\nm)^*\text{ and }
S_{\np}^{\n}:=U(\n)\ten_{U(\np)}U(\np)^*.
$$
Then the $\g$-modules $S_{\nm}$ and $U(\g)\ten_{U(\n)}S_{\nm}^{\n}$
(resp.
the $\g$-modules $S_{\np}$ and $U(\g)\ten_{U(\n)}S_{\np}^{\n}$)
are naturally isomorphic.

The induction functor $\Ind_{\nm}^{\n}:\ \nm\mod\map\n\mod$
takes $U(\nm)^*$ to $S_{\nm}^{\n}$ and provides an inclusion of algebras
$U(\n)\subset\End_{\n}(S_{\nm}^{\n}$):
$$
n\cdot a\ten f:=a\ten n\cdot f,\text{ where }n\in\nm, f\in U(\nm)^*,a\in U(\n),
\ n\cdot f(n'):=f(n'n), \n'\in\nm.
$$
\sssn
Consider the following action of $\np$ on $U(\n)\ten U(\nm)^*$:
$$
n^+\cdot a\ten f:=an^+\ten f+a\ten[n^+,f],\text{ where }
n^+\in\np,f\in U(\nm)^*,a\in U(\n),\ [n^+,f](u):=f([u,n^+]),\ u\in U(\nm).
$$
\sssn
\Lemma
The defined action of $\np$  on $U(\n)\ten U(\nm)^*$
commutes with the left regular action of $U(\n)$.
Moreover it is well defined on $S_{\nm}^{\n}$. Along with the action of
$\nm$ on $S_{\nm}^{\n}$ it defines the inclusion of algebras $U(\n)\subset
\End_{\n}(S_{\nm}^{\n})$.

\dok
The first statement is obvious. Let $a\in U(\n)$, $f\in U(\nm)^*$,
$n^+\in\np$, $n^-\in\nm$. Note that $[n^-,n^+]\in\nm$. Thus we have
\begin{gather*}
n^+\cdot an^-\ten f=
an^-n^+\ten f+a^n-\ten[n^+,f]=
an^-n^+\ten f+a\ten n^-[n^+,f]\\=
an^+n^-\ten f+a\ten [n^-,n^+]f-a\ten [n^-,n^+]f+a\ten [n^+,n^-f]=
n^+\cdot a\ten n^-f.
\end{gather*}
To prove the third statement note that
\begin{gather*}
n^+\cdot n^-\cdot a\ten f=
an^+\ten n^-f+a\ten[n^+,n^-f]\\=
an^+\ten n^-f+a\ten n^-[n^+,f]+a\ten[n^+,n^-]f=
n^-\cdot n^+\cdot a\ten f+[n^+,n^-]\cdot f.\qed
\end{gather*}
The following statement is crusial in the proof of the Theorem.

\sssn
\Lemma
Thre exists an isomorphism of left $\n$-modules $U(\n)^*\cong S_{\nm}^{\n}
\ten_{U(\np)}U(\np)^*$

\dok
Let $f^+\in U(\np)^*$, $f^-\in U(\nm)^*$. Denote the element
$(1\ten f^-)\ten f^+\in S_{\nm}^{\n}\ten_{U(\np)}U(\np)^*$ by $f^-\ten
f^+$.  Then such elements form a base of the vector space
$S_{\nm}^{\n}\ten_{U(\np)}U(\np)^*\cong U(\nm)^*\ten U(\np)^*$. We calculate
the action of $\n$ on $S_{\nm}^{\n}\ten_{U(\np)}U(\np)^*$ in this base.  For
$n^-\in\nm$ we have
$$
n^-(f^-\ten f^+)=(n^-f^-)\ten f^+.
$$
For $n^+\in\np$ we have
$$
n^+(f^-\ten f^+)=(n^+\ten f^-)\ten f^+=(n^+\cdot(1\ten f^-))\ten f^+
-[n^+,f^-]\ten f^+=
f^-\ten n^+f^+-[n^+,f^-]\ten f^+.
$$
Now recall that $\nm$ is an ideal in $\n$ and the actions  of the subalgebras
$\nm$ and $\np$ on $U(\n)^*\cong U(\nm)^*\ten U(\np)^*$ are given by these
very formulas.  \qed

We return to the situation of (\ref{main2}). Note that each module
$S_{\n^k}$ is in fact a $\g$-bimodule by Lemma~\ref{act}.

\sssn
\Lemma  \label{iterate}
The left $\g$-modules $S_{\n}$ and $S_{\n^0}\ten_{U(\g)}\ldots\ten_{U(\g)}
S_{\n^{\operatorname{top}-1}}$ are naturally isomorphic to each other.

\dok
We prove by induction by $\operatorname{top}-k$ that
the $\g$-modules $S_{F^k\n}$ and $S_{F^{k+1}\n}\ten_{U(\g)}S_{\n^k}$
are naturally isomorphic to each other.
Note that each time the induction hypothesis provides the $\g$-bimodule
structure on $S_{F^{k+1}\n}$:
$$
U(\g)\hookrightarrow\End_{\g}(S_{F^{\operatorname{top}-1}\n})
{\overset{\cdot\ten_{U(\g)}S_{\n^{\operatorname{top}-2}}}
{\map}}
\ldots
{\overset{\cdot\ten_{U(\g)}S_{\n^{k+1}}}
{\map}}
\End_{\g}(S_{F^{k+1}\n}).
$$
Thus by the previous Lemma
$$
S_{F^{k+1}\n}\ten_{U(\g)}S_{\n^k}\cong
U(\g)\ten_{U(F^k\n)}S_{F^{k+1}\n}^{F^k\n}\ten_{U(\n^k)}U(\n^k)^*
\cong
U(\g)\ten_{U(F^k\n)}U(F^k\n)^*=S_{F^k\n}.
$$
The Lemma is proved.\qed

The statement of the Theorem follows immediately from the previous Lemma.\qed

\section{Affine root systems and affine Weyl groups}

\subsection{General setting}
In this section we recall some basic facts about affine Weyl groups.
In our notations we follow mostly
[L1], 1.1 -- 1.5. A detailed exposition of the subject can be found in [K],
Chapter 6.

\sssn
Let
$(a_{ij})_{i,j\in\afr}$ be an irreducible (untwisted) affine Cartan matrix.
There are uniquely defined strictly positive integers $d_i,r_i,r'_i,\ i\in\afr$
such that

(i) $d_ia_{ij}=d_ja_{ji}$ for all $i,j$ and $d_i=1$ for some $i$,\\
(ii) $\sum_ir_ia_{ij}=0$ for all $j$ and $r_0=1$,\\
(iii) $\sum_jr'_ja_{ij}=0$ for all $i$ and $r'_0=1$.

\sssn   \label{finite}
In particular $(a_{ij})_{i,j\in\underline{r}}$ is a Cartan matrix of the finite
type.

Let $D=\max_id_i$, we have $D\in\{1,2,3\}$ and for each $i$,  $d_i$ is equal
either to $1$ or $D$. We define $\hat{d}_i$ by $d_i\hat{d}_i=D$ for  all
$i\in\afr$.

\sssn
Let $V$ be a $\R$-vector space with basis $\{h_i| i\in\afr\}$ and let $V'$ be
the dual vector space. We denote by $\langle\cdot,\cdot\rangle;\
V\times V'\map\R$ the canonical bilinear pairing. Let $\{\omega_i|i\in\afr\}$
be the basis of $V'$ dual to $\{h_i|i\in\afr\}$. Define $\{\alpha_j|
j\in\afr\}$ by $\langle h_i,\alpha_j\rangle=a_{ij}$. Consider a vector
$c:=\sum_i r_ih_i\in V$.  Then we have $\langle c,\alpha_j\rangle=0$ for all
$j\in\afr$ and $\sum_ir'_i\alpha_i=0$.

\sssn
For $i\in\afr$ we define reflections
$$
s_i:\ V\map V,\ s_i(y):=y-\langle y,\alpha_i\rangle h_i,\text{ and }
s_i:\ V'\map V',\ s_i(x):=x-\langle h_i,x\rangle \alpha_i.
$$
The subgroup $W\subset GL(V)$ generated by the reflections
$s_i,\ i\in\afr$, is called the affine Weyl group. We identify $W$
with the subgroup
in $GL(V')$ generated by $s_i,\ i\in\afr$. The subgroup
$\overline{W}\subset W$ generated by the reflections $s_i,\
i\in\underline{r}$, is called the (finite) Weyl group corresponding to the
Cartan matrix $(a_{ij})_{i,j\in\underline{r}}$.

\sssn
The set $R$ (resp. $\overline{R}$) of vectors of $V$ of the form
$w(h_i)$ for some $i\in\afr$ and $w\in W$
(resp. of the form
$w(h_i)$ for some $i\in\underline{r}$ and $w\in\overline{W}$)
is called the affine root system (resp. the finite root system).
Let $R'$ (resp. $\overline{R}'$) be the set of vectors of $V'$ of the form
$\omega(h_i)$ for some $i\in\afr$ (resp. for some $i\in\overline{r}$) and
some $\omega\in \overline{W}$). The assignment $h_i\mapsto \alpha_i$ extends
uniquely to a map
$$
h\mapsto h',\ R\map R'\text{, such that  } \omega(h)'=\omega(h'),\
\omega\in W,\ h\in V.
$$
The map restricts to bijection of $\overline{R}$ to  $\overline{R}'$.

There is a unique function $h\mapsto d_h$ on $R$ such that it is
$W$-invariant and $d_{h_i}=d_i$ for all $i\in\afr$. We define $\hat{d}_h$ by
$d_h\hat{d}_h=D$ for all $h\in R$.

\sssn
\Rem

(i) $R=\{\alpha+\hat{d}_{\alpha}mc|\alpha\in\overline{R},m\in\Z\}$;\\
(ii) $R'=\overline{R}'$;\\
(iii) $(\alpha+\hat{d}_{\alpha}mc)'=\alpha'$ for all $\alpha\in
\overline{R},\ m\in\Z$.

For $h\in R$ we denote by $s_h$ the element of $W$  given by the reflection
in $V$ (resp. in $V'$)
$$
s_h(y)=y-\langle y,h'\rangle h\text{ (resp. }s_h(x)=x-\langle h,x\rangle h').
$$
For any $\alpha\in\overline{R}$ and $m\in\Z$ we set $s_{\alpha,m}=s_h\in W$,
where $h=\alpha+\hat{d}_{\alpha}mc$.

Consider the weight lattice $P\subset V'$ generated by the set
$\{\omega_i|i\in\afr\}$.
Let
$Q'\subset V'$ be a free abelian group
generated by the set
$\{\alpha_i|i\in\underline{r}\}$ and let
$Q''\subset Q'$ be
a free abelian group
generated by the set
$\{\hat{d}_i\alpha_i|i\in \underline{r}\}$. For $z\in Q''$ consider  a
transvection $\theta_z:\ V'\map V'$ given by $\theta_z(x)=x+\langle
c,x\rangle z$.

\sssn
\Lemma
For $\alpha\in\overline{R}$ and $m\in\Z$ we have $s_{\alpha,0}\circ
s_{\alpha,m} =\theta_{\hat{d}_{\alpha}m\alpha'}$.  \qed

In particular $\theta_z\in W$ for any $z\in Q''$. Consider the map of the sets
$
\theta:\ Q''\map W,\
z\mapsto \theta_z,
$
and denote its image by $T\subset W$.

\sssn
\Lemma

\qquad(i) The map $\theta$ is an injective homomorphism of groups;

\qquad(ii) $T$ is a normal subgroup in $W$;

\qquad(iii) $W$ is a semidirect product of $T$ and $\overline{W}$.
\qed

Let $V^+$ (resp. $V^-$)be the set of all vectors in $V$  such that all their
coefficients with respect to the basis $\{h_i|i\in \afr\}$ are nonnegative
(resp. nonpositive).  Set $R^+=R\cap V^+$, $R^-=R\ V^-$, $\overline{R}\cap
V^+$.  We have
$R^+=\{\alpha+\hat{d}_{\alpha}mc|\alpha\in\overline{R},m>0\}\sqcup
\overline{R}^+$.

The height of a positive root $\alpha=\underset{i\in\afr}{\sum}b_i\alpha_i$
is defined as follows: $\hgt\alpha=\underset{i\in\afr}{\sum}b_i$.

\sssn
As usual define the weight $\rho\in V'$ by
$\rho(h_i)=1$ for all $i\in\afr$. Consider the {\em dot} action of
$W$ on $V'$:
$$
w\cdot\lambda=w(\lambda+\rho)-\rho,\ w\in W,\ \lambda\in V'.
$$
Then  the dot action of the Weyl group preserves the sets
$P_k:=\{\lambda\in P|\langle c,\lambda\rangle=k\}$.

From now on we denote the set of dominant weights at the level $k$
by $P^+_k$:
$$
P^+_k:=\{\lambda\in P_k|\langle h_i,\lambda\rangle\ge0,\ i\in\afr\}.
$$
\subsection{Semi-infinite Bruhat order on the Weyl group.}
Recall that  the length of an element of the affine Weyl group $w$ is defined
as follows:
$$
\lth(w):=\sharp\{\alpha\in R^+|w^{-1}(\alpha)\in R^-\}.
$$
\Rem
The length of $w\in W$ is equal also to the minimal possible length of
expression of $w$ via the generators $s_i,\ i\in\afr$ (the length of a
reduced expression).

\sssn            \label{length}
For $w\in W$ consider a finite set
$R_w:=\{\alpha\in R^-|w(\alpha
)\in R^+\}$. Let $w_1$ and $w_2$ be elements of the Weyl group
such that  $\lth(w_1)+\lth(w_2)=\lth(w_2w_1)$. Then $R_{w_2w_1}=
R_{w_2}\sqcup R_{w_1}$. Thus $R_{w_2}\cap-w_1(R_{w_1})=\emptyset$, where the
set $-S\subset R$ consists of elements opposite to the ones of $S\subset R$.

Recall that we have
$\underset{\alpha\in R_w}{\sum}\alpha=\rho-w^{-1}(\rho)$.

\sssn
We say that $w'$ follows $w$ in the Weyl group if there exist a reduced
expression of $w'$ and  $p\in\{1,\ldots,\lth w'\}$ such that
$$
w'=s_{i_1}\ldots s_{i_{\lth w'}},\ w=s_{i_1}\ldots s_{i_{p-1}}s_{i_{p+1}}\ldots
s_{i_{\lth w'}}
$$
and $\lth w=\lth w'-1$.

Recall that the usual Bruhat order on the Weyl group is the partial order on
$W$ generated by the relation ``$w'$ follows $w$ in $W$''. It is denoted by
$\ge$.

\sssn
\Lemma
The relation
$$
\{\text{there exists }\lambda_0\in Q^{\prime\prime+}\text{ such that for any }
\lambda\in Q^{\prime\prime+}\text{ we have }
\theta_{\lambda}\theta_{\lambda_0}w'
\ge
\theta_{\lambda}\theta_{\lambda_0}w\}
$$
is a partial order on $W$.
\qed

\sssn
\Def
We call this partial order the semi-infinite Bruhat order on $W$ and denote it
by $\ge^{\si}$.

\sssn
\Rem
The semi-infinite Bruhat order defined above in fact coincides with the
partial order on the affine Weyl group defined in [L2], Section 3, in terms
of combinatorics of {\em alcoves} in $\overh^*$. However we will not need the
comparison statement. Further details on the partial order can be found in
[L2].

\sssn  \label{tlength}
Denote the set $\{\alpha\in R^+|\alpha=\beta+\hat{d}_{\beta}mc,\
\beta\in\overr^+,m\ge 0\}$ (resp.  the set $\{\alpha\in
R^+|\alpha=\beta+\hat{d}_{\beta}mc,\ \beta\in\overr^-,m> 0\}$) by $R^{\frac
\infty 2+}$ (resp.  by $R^{\frac \infty 2-}$).  Following [FF] we introduce
the semi-infinite length function on the affine Weyl group as follows:
$$
\lth^{\si}(w):=
\sharp\{\alpha\in R^{\si+}|w(\alpha)\in R^-\}
-
\sharp\{\alpha\in R^{\si-}|w(\alpha)\in R^-\}.
$$
Next we consider  the twisted length function on the affine Weyl group.
For $u,w\in W$ set
$$
\lth^w(u):=\lth(w^{-1}u)-\lth(w^{-1}).
$$
In particular for $\mu\in -Q^{\prime\prime+}$ and $\nu\in Q^{\prime\prime+}$
we have $\lth^{\theta_{\mu}}(\theta_{\nu})=\lth(\theta_{\nu})$.

\sssn
\Lemma
\label{maincomb} For every $w_1,w_2\in W$ there exists $\mu_0\in
-Q^{\prime\prime+}$ such that for every $\mu\in -Q^{\prime\prime+}$ we have
$$
\lth(\theta_{-\mu}\theta_{-\mu_0}w_1)-
\lth(\theta_{-\mu}\theta_{-\mu_0}w_2)=
\lth^{\si}(w_1)-
\lth^{\si}(w_2).
$$
In particular for every $\mu\in -Q^{\prime\prime+}$ we have
$\lth^{\theta_{\mu}\theta_{\mu_0}}(w)=\lth^{\si}(w)$.
\qed

\section{Twisted Verma modules over affine Lie algebras and twisting functors}

\subsection{Affine Lie algebras.}
Fix a Cartan matrix $(a_{ij})_{i,j\in\underr}$ like in~\ref{finite} and
consider the corresponding semisimple Lie algebra $\overg$. In particular
$\overr$ is the root system of $\overg$.  Thus $\overg=\overline{\g}^+\oplus
\overh\oplus \overg^-$, where $\overh$ denotes the Cartan subelgebra in
$\overg$, and
$$
\overline{\g}^+=\underset{\alpha\in\overline{R}^+}
{\bigoplus}\overline{\g}_{\alpha},\
\overline{\g}^-=\underset{\alpha\in\overline{R}^-}
{\bigoplus}\overline{\g}_{\alpha}.
$$
\sssn
Recall that the affine Lie algebra $\g$ for the semisimple Lie algebra
$\overg$ is defined as a central extension of a loop algebra ${\cal
L}\overg:=\overg\ten\CC$ $[t,t^{-1}]$.  Namely,
$\g:={\cal L}\overg\oplus\CC K$, and the bracket in $\g$ is defined as
follows:
$$
 [g_1\ten t^n,g_2\ten t^m]=[g_1,g_2]\ten
 t^{n+m}+\delta_{n+m,0}nB(g_1,g_2)K,
$$
where
$g_1,g_2\in\overg,n,m\in \Z$, and $B(\cdot,\cdot)$ denotes the
Killing form of $\overg$.

\sssn
It is well known that the affine Lie algebra $\g$ is a Kac-Moody Lie
algebra with the Cartan matrix $(a_{ij})_{i,j\in \{0\}\cup\underr}$ and
the root system $R$.
Thus
$\g=\n^+\oplus\h\oplus \n^-$, where $\h:=\overh\oplus\CC K=V\ten_{\R}\CC$
and $$ \n^+=\underset{\alpha\in
R^+}{\bigoplus}\g_{\alpha}=\overline{\g}^+\oplus \overline{\g}\ten
t\CC[t],\ \n^-=\underset{\alpha\in
R^-}{\bigoplus}\g_{\alpha}=\overline{\g}^-\oplus \overline{\g}\ten
t^{-1}\CC[t^{-1}].
$$
\sssn
Recall the definition of the Chevalley generators for $\g$. Set
$e_{\alpha_i}:=e_i$, $e_{-\alpha_i}:=f_i$,
$i=1,\ldots,r$, where $e_i,f_i$
are the Chevalley generators for $\overline{\g}$,
$[e_i,f_i]=h_i$. Let $\alpha_{\operatorname{top}}$
be the highest root of the root system $\overr$.
Choose a vector $h_{\operatorname{top}}$ in
$[\overline{\g}_{\alpha_{\operatorname{top}}},
\overline{\g}_{-\alpha_{\operatorname{top}}}]\subset\overh$ such that
$\alpha_{\operatorname{top}}(h_{\operatorname{top}})=2$.
Fix
$e_{\operatorname{top}}\in\overline{\g}_{\alpha_{\operatorname{top}}}$
and
$f_{\operatorname{top}}\in\overline{\g}_{-\alpha_{\operatorname{top}}}$
such that
$[e_{\operatorname{top}},f_{\operatorname{top}}]=h_{\operatorname{top}}$.
Then set
$$
e_{\alpha_0}:=f_{\operatorname{top}}\ten t,\
e_{-\alpha_0}:=e_{\operatorname{top}}\ten t^{-1}.
$$
\sssn \label{gr}
We introduce a grading on $\g$ putting $\deg e_i=1$,
$\deg f_i=-1$ and $\deg h_i=0 $ for all $i\in\afr$.

\subsection{Categories of $\g$-modules.}
Consider the category $\M$
of $\h^*\times\Z$-graded $\g$-modules $M=\underset{\lambda\in\h^*,t\in\Z}
{\bigoplus}M_{\lambda,t}$ such that

(i) for every $i\in\afr$ the   Chevalley generators
$
e_i:\ M_{\lambda,t}\map M_{\lambda+\alpha_i,t+1},\
f_i:\ M_{\lambda,t}\map M_{\lambda-\alpha_i,t-1};
$\\
(ii) every $h\in\h$ acts on $M_{\lambda,t}$ by scalar $\lambda(h)=
\langle h,\lambda\rangle$.

Morphisms in $\M$ are  morphisms of $\g$-modules that preserve
$\h^*\times\Z$-gradings.

Fix a nonnegative integer $k\in\Z_{\ge 0}$. Denote by $U_k(\g)$ the quotient
algebra of the universal enveloping algebra of $\g$ by the relation $K-k=0$.
Consider a full subcategory $\M_k$ in $\M$ of modules $M$ such that $K$ acts
by the scalar $k$ on $M$. Denote by $\supp M$ the set
$\{\lambda\in\h^*|M_{\lambda,t}\ne0\text { for some }t\in\Z\}$. Then for
every $M\in\M_k$ the set $\supp M$ is contained in the set
$\h^*_k:=\{\lambda\in\h^*|\langle c,\lambda\rangle=k\}$.  The set $\h^*_k$ is
preserved both by the linear action of the affine Weyl group on $\h^*$ and by
the dot action.

\sssn
Let $D$ be the contragradient duality functor:
$$
D:\ \M\map\M,\ D(M):=\underset{\lambda\in\h^*,t\in\Z}{\bigoplus}
M_{-\lambda,-t}^*,\ M\in\M.
$$
The left action of $\g$ on $D(M)$ is provided by the composition of the usual
right action of $\g$ on the dual module with the antipode and the Chevalley
involution of $\g$.  In particular $D$ preserves $\M_k$.

\sssn
We define the character of a $\g$-module $M\in\M$ such that $\dim
M_{\lambda,t}<\infty$ for all $\lambda\in\h^*,t\in\Z$ by
$$
\ch M:=\sum_{\lambda\in\h^*,t\in\Z}\dim M_{\lambda,t}e^{\lambda}q^t.
$$
Here $q$ is a formal variable and $e^{\lambda}$ is a formal expression.

\sssn
As usual let $\O$ denote the category of $\h^*$-graded
$\g$-modules $M=\underset{\lambda\in\h^*}{\bigoplus}M_{\lambda}$ such that

(i) for every $i\in\afr$ the   Chevalley generators
$
e_i:\ M_{\lambda}\map M_{\lambda+\alpha_i,},\
f_i:\ M_{\lambda}\map M_{\lambda-\alpha_i,};
$\\
(ii) every $h\in\h$ acts on $M_{\lambda}$ by scalar $\lambda(h)=
\langle h,\lambda\rangle$;\\
(iii) $\dim M_{\lambda}<\infty$ for all $\lambda\in\h^*$;\\
(iv) there exist $\lambda_1,\ldots,\lambda_n\in\h^*$ such that
$$
M_{\mu}\ne0\text{ only when }\mu\in\lambda_1+R^-\cup\ldots\cup\lambda_n+R^-.
$$
\sssn
We recall the definition of the Casimir endomorphism $\Gamma_M:\ M\map M$
for $M\in\O$. Choose a homogenious base
$\{e_{\alpha}|\alpha\in R\}\cup\{e_{0,i}|i\in\afr\}$
in $\g$ and let
$\{e^{\alpha}|\alpha\in R\}\cup\{e_0^i|i\in\afr\}$ be the dual base with
respect to the Killing form on $\g$.  For $M\in\O$ we  define the Casimir
operator by
$$
\Gamma_M|_{M_\nu}=(\nu+\rho,\nu+\rho)\operatorname{Id}_{M_{\nu}}
+2\sum_{\alpha\in R^+}e^{\alpha}e_{\alpha}
$$
where $e^{\alpha}e_{\alpha}\in U(\g)$.
It is known that $\Gamma_M$ commutes with the action of $\g$.
For $M\in\O$ and $\kappa\in\CC$ set
$$
M_{\kappa}:=\{m\in M|\text{ there exists }n\in\N\text{ such that
}(\Gamma_M-\kappa)^nm=0\}
$$
\sssn
\Lemma
$M_{\kappa}$ is a submodule in $M$ and
$M=\underset{\kappa\in\CC}{\bigoplus}M_{\kappa}$.  \qed

\subsection{Affine BGG resolution.}
Let $\b^+$ (resp. $\b^-$) be the positive (resp. negative) Borel subalgebra
in $\g$,
$\b^+:=\h\oplus\underset{\alpha\in R^+}{\bigoplus}\g_{\alpha}$
(resp.
$\b^-:=\h\oplus\underset{\alpha\in R^-}{\bigoplus}\g_{\alpha}$).
Then $\n^+$ is an ideal in $\b^+$ and $\b^+/\n^+=\h$. For the one
dimensional $\h$-module $\CC(\lambda)$ the  Verma module $U_k(\g)\ten_{U(\b^+)}
\CC(\lambda)$ is denoted by $M(\lambda)$.
Evidently $M(\lambda)$ is free as a $U(\n^-)$-module and
belongs to $\cal O$ and to $\M_k$ if we put the
$\Z$-grading of the highest weight
vector $v_{\lambda}\in M(\lambda)$ equal to zero.

Recall the construction of the BGG resolution
in the case of affine Lie algebras (see e.~g [RW]).
For a fixed positive integral level $k$  choose a dominant weight
$\lambda\in P_k^+$.
Consider the unique simple quotient module of $M(\lambda)$ denoted by
$L(\lambda)$. Then there
exists a left resolution of $L(\lambda)$ that consists of direct sums of  Verma
modules. Namely recall  first the construction of
the standard resolution of the trivial $\g$-module
$\underc$ relative to $\b^+$.

\subsubsection{Standard resolutions.} \label{stres}
Consider the standard complex $U(\g)\ten \Lambda\bul(\g)\ten\underc$
for the computation of the Lie algebra homology of
$\g$ with coefficients in $U(\g)$. Clearly
$$
H^{\ne0}(U(\g)\ten\Lambda\bul(\g)\ten\underc)=0,\
H^0(U(\g)\ten\Lambda\bul(\g)\ten\underc)=\underc
$$
and $
U(\g)\ten\Lambda\bul(\g)\cong
U(\g)\ten\Lambda\bul(\b)\ten\Lambda\bul(\g/\b)$ as a vector space.
We introduce a filtration of the complex as follows:
$$
F^m(U(\g)\ten\Lambda^n(\g)):=\underset{{p+q=n}\atop{q\le m}}{\bigoplus}
U(\g)\ten\Lambda^p(\b)\ten\Lambda^q(\g/\b).
$$
Clearly the differential in the complex preserves the filtration. Consider
the corresponding spectral sequence. We have
$$
E_1^{p,q}=H_{-p}(\b,U(\g)\ten\Lambda^{-q}(\g/\b))=
\delta_{p,0}U(\g)\ten_{U(\b)}\Lambda^{-q}(\g/b).
$$
Thus the spectral sequence converges in $E_2$. On the other hand
$E_{\infty}^{p,q}=\delta_{p,0}\delta_{q,0}\underc$ and we obtain the required
standard resolution of the trivial $\g$-module $\underc$ with relative to
$\b$:
\begin{gather*}
D\bul:=U(\g)\ten_{U(\b^+)}\Lambda\bul(\g/\b^+),\ d:\ D^{-k}\map D^{-k+1},\\
d(u\ten
\overline{a}_1\wedge\ldots\wedge\overline{a}_k)
=\sum_{i=1}^k(-1)^i
ua_i\otimes\overline{a}_1\wedge\ldots\wedge\overline{a}_{i-1}
\wedge\overline{a}_{i+1}\wedge\ldots\wedge\overline{a}_k\\
+\sum_{i<j}(-1)^{i+j}u\otimes\overline{[a_i,a_j]}\wedge\overline{a}_1\wedge
\ldots\wedge\overline{a}_{i-1}\wedge\overline{a}_{i+1}\wedge\ldots\wedge
\overline{a}_{j-1}\wedge\overline{a}_{i+1}
\wedge\ldots\wedge\overline{a}_k,
\end{gather*}
where $\overline{a}$ denotes the image of $a$ in $\g/\b^+$.

\sssn
\Rem
Suppose we have two subalgebras $\b_1\subset\b_2\subset\g$. Then the identity
map $U(\g)\ten\Lambda\bul(\g)\map U(\g)\ten\Lambda\bul(\g)$ is a morphism of
filtered complexes:
$
F^m_{\b_1}(U(\g)\ten\Lambda\bul(\g))\subset
F^m_{\b_2}(U(\g)\ten\Lambda\bul(\g))
$.
Thus we have a canonical morphism of the standard resolutions
$D\bul_{\b_1}\map D\bul_{\b_2}$.

\sssn
Consider the standard resolution of the simple module $L(\lambda)$ relative
to $\b^+$:
$
D\bul(L(\lambda)):=D\bul\ten L(\lambda)
$
with the structure of $\g$-module using comultiplication on $U(\g)$.
Then
$D\bul(L(\lambda))=\underset{\kappa\in\CC}{\bigoplus}
\left(D\bul(L(\lambda))\right)_{\kappa}$.

Taking the direct summand with $\kappa=(\lambda+\rho,\lambda+\rho)$ we obtain
a resolution of $L(\lambda)$ filtered by Verma modules. The folowing statement
says that the filtration splits.

\sssn
\Theorem
(see [RW] 9.7)
There exists a resolution
$B\bul(\lambda)=\left(D\bul(L(\lambda))\right)_{(\lambda+\rho,\lambda+\rho)}$
of $L(\lambda)$ in $\O$ and in $\M_k$  of the form
$$
\ldots
\map
\underset{{w\in W,}\atop{\lth(w)=m}}{\bigoplus}
M(w\cdot\lambda)\langle-\hgt(\lambda-w\cdot\lambda)\rangle
\map\ldots\map
\underset{{w\in W,}\atop{\lth(w)=1}}{\bigoplus}
M(w\cdot\lambda)\langle-\hgt(\lambda-w\cdot\lambda)\rangle
\map
M(\lambda)\map L(\lambda)\map 0.
$$
Here as usual $\langle\cdot\rangle$ denotes the shift of the $\Z$-grading
in $\M$.
\qed

\subsection{Twisting functors}
Recall that the affine Weyl group has a geometric description in terms of
affine Lie groups. Let $\hat{{\cal L}G}$ be the central extension of the
loop group of $\overg$. Let $T\subset \hat{{\cal L}G}$ be the Cartan subgroup.
Denote by $N(T)$ its normalizer in $\hat{{\cal L}G}$.

\sssn
\Lemma
The affine Weyl group $W$ is isomorphic to $N(T)/T$.
\qed

For every $w\in W$ fix its representative $\dot w\in N(T)$. The normalizer
of $T$ acts
on $\g$
by adjunction
and this action commutes  with the
bracket in $\g$ and shifts the weight decomposition.
In particular we have maps
$$
\dot w: \g\map\g, \
\g_{\alpha}\map\g_{w(\alpha)}.
$$
\sssn       \label{tf}
We introduce the functors of the
twist of the $\h^*$-grading.  We define $T_w:\ \M_k\map\M_k$ as follows.
For $M\in\M_k$ set
\begin{gather*}
T_w(M)=\underset{\lambda\in\h^*,t\in\Z}{\bigoplus}T_w(M)_{\lambda,t},\
T_w(M)_{\lambda,t}:=M_{w(\lambda),t+\hgt(w(\lambda)-\lambda)},\
e_{\alpha}\in\g_{\alpha},\ m_{\lambda}\in T_w(M)_{\lambda,t}\text{ then }\\
e_{\alpha}\cdot m_{\lambda}:=\dot w(e_{\alpha})(m_{\lambda})\in
M_{w(\lambda+\alpha),t+\hgt (w(\lambda+\alpha)-\lambda)}
=M_{\lambda+\alpha,t+\hgt (w(\lambda+\alpha)-\lambda-\alpha)+\hgt\alpha}
=T_w(M)_{\lambda+\alpha,t+\hgt\alpha}.
\end{gather*}
Evidently $T_w$ is an equivalence of categories and the opposite functor is
$T_{w^{-1}}$.

\sssn
Fix $w\in W$. Let $\n_w$ be  a  finite
dimensional nilpotent subalgebra in $\g$:
$$
\n_w:=\n^-\cap\dot w^{-1}(\n^+)=\underset{\alpha\in R_w}{\bigoplus}\g_{\alpha}.
$$
Consider the semiregular $\g$-module with respect to $\n_w$ at the level $k$:
$$
S_{\n_w}=\Ind_{U(\n_w)}^{U_k(\g)}\Coind_{\CC}^{U(\n_w)}\underc.
$$
Then by \ref{main2} there exists an inclusion of algebras $\sigma_w:\
U_k(\g)\hookrightarrow \End_{\g}(S_{\n_w})$.

Like in the third section, for $\alpha\in R_w$ consider the $\g$-module
$S_{\n_{\alpha}}:=U_k(\g)\ten_{U(\n_{\alpha})}U(\n_{\alpha})^*$.
By~\ref{main2} there also exists an inclusion of algebras $\sigma_{\alpha}
:\ U_k(\g)\hookrightarrow\End_{\g}(S_{\n_{\alpha}})$ (see~\ref{exp}).
As before identify $U(\n_{\alpha})^*$ with $\CC[x_{\alpha}]$ so that
$e_{\alpha}\in\n_{\alpha}$ acts by $\pl/\pl x_{\alpha}$.

\sssn
\Lemma                                             \label{h}
For every $h\in\h$ we have
$\sigma_{\alpha}(h)=h-\alpha(h)e_{\alpha}\ten x_{\alpha}$.

\dok
Follows immediately from~\ref{act}.
\qed

\sssn
Consider a functor
$$
S_w:\ \M_k\map U_k(\g)\mod,\ S_w(M):=S_{\n_w}\ten_{U_k(\g)}M.
$$
We will show that in fact $S_w$ can be dafined as a functor from $\M_k$
to $\M_k$.

\Rem                             \label{1}
Note that $S_w(M)\cong U(\n_w)^*\ten_{U(\n_w)}M$ as a vector space (and even
as a $U(\n_w)$-module).

Let $f\in U(\n_w)^*_{-\lambda}=(U(\n_w)_{\lambda})^*$ and
$m\in M_{\mu,t}$. Then set the $\h^*\times\Z$-grading of $f\ten m\in S_w(M)$
equal to $(\mu-\lambda+\rho-w^{-1}(\rho),t+\hgt(\rho-w^{-1}(\rho)))$.  The
grading is well defined with respect to the  action of  $\g$. It remains to
check that $S_w(M)$ is $\h$-semisimple. But this follows from~\ref{h}
and~\ref{iterate}.  We have proved the following statement.

\sssn
\Lemma
$S_w$ defines a functor $\M_k\map\M_k$.
\qed

\sssn
\Def
For $w\in W$ the functor of twist by $w$
$$
\Phi_w:=T_w\circ S_w,\ \Phi_w:\ \M_k\map\M_k.
$$
Let us describe the image of a  Verma module under $\Phi_w$.

\sssn
\Lemma            \label{weight}
$\ch \Phi_w(M(\lambda))=\ch M(w\cdot\lambda)\langle
-\hgt(\lambda-w\cdot\lambda)\rangle$.

\dok
Follows immediately from the definition of the $\h^*\times\Z$-grading on
$T_w\circ S_w(M(\lambda))$.
\qed

In particular the highest weight vector $1\ten
v_{\lambda}\in\Phi_w(M(\lambda))$ has the weight $w\cdot\lambda$.

\sssn
\Def
We call the $\g$-module $M_w(w\cdot\lambda):=\Phi_w(M(\lambda))\langle
\hgt(\lambda-w\cdot\lambda)\rangle$ the {\em twisted Verma module} of the
weight $w\cdot\lambda$.

\sssn
\Rem  \label{semireg}
Consider the nilpotent subalgebras
$\dot w(\n_w)=\dot w(\n^-)\cap\n^+\subset
\dot w(\n^-)\subset\g$.
Then when restricted to $\dot w(\n^-)$ the  twisted Verma module
$M_w(w\cdot\lambda)$ is isomorphic to the $\dot w(\n^-)$-module
$S_{\dot w(\n_w)}^{\dot w(\n^-)}$.

\subsection{Twisted BGG resolutions.}
Fix $w\in W$. Consider the complex in $\M_k$ of the form
$\Phi_w(B\bul(\lambda))$. Using   the fact that
$B\bul(\lambda)$ is a $U(\n_w)$-free resolution of $L(\lambda)$
we see that up to a $\Z$-grading shift
$$
H^{-\bullet}
(\Phi_w(B\bul(\lambda)))=\Tor^{\n_w}(U(\n_w)^*,L(\lambda)).
$$
as vector spaces.

\sssn
\Prop      \label{phi}
$H^{\ne -l(w)}
(\Phi_w(B\bul(\lambda)))=0,\
H^{ -l(w)}
(\Phi_w(B\bul(\lambda)))=L(\lambda)$ as a $\g$-module.

\subsection{Proof of the Proposition~\ref{phi}.}
We present first several Lemmas. Let $\n$ be a negatively graded finite
dimensional Lie algebra.

\sssn
\Lemma
$\tor_{\ne\dim\n}^{\n}(U(\n)^*,\underc)=0,\
\tor_{\dim\n}^{\n}(U(\n)^*,\underc)=\CC$.

\dok
Consider the standard complex for the computation of $\n$-homology
$K\bul(U(\n)^*):=U(\n)^*\ten\Lambda\bul(\n)\ten\underc$.  Then the PBW
filtration on $U(\n)$ induces filtrations both on $U(\n)^*$ and on
$\Lambda(\n)$ (such that for any base vector $e_i$ of $\n\subset\Lambda(\n)$
we have $e_i\in F^1\Lambda(\n)$ and for any dual base vector $e_i^*\in
U(\n)^*$ we have $e_i^*\in F^{-1} U(\n)^*$). The differential in
$K\bul(U(\n)^*)$ preserves the filtration. The complex $\gr^FK\bul(U(\n)^*)$
is isomorphic to the coKoszul complex  $S(\n)^*\ten\Lambda\bul(\n)$  of the
algebra $\Lambda(\n)$.  Thus $$
H^{\ne-\dim\n}(\gr^FK\bul(U(\n)^*))=0, H^{-\dim\n}(\gr^FK\bul(U(\n)^*))=\CC.
$$ Now the grading on $\n$ induces a grading on $K\bul(U(\n)^*)$.  The
spectral sequence of the filtration $F$ is finite in each grading component
of the complex $K\bul(U(\n)^*)$. Thus it converges.  \qed

\sssn
\Cor
Let $M$ be  a finite dimensional graded $\n$-module. Then
$$
\tor_{\ne\dim\n}^{\n}(U(\n)^*,M)=0,\
\tor_{\dim\n}^{\n}(U(\n)^*,M)=M
$$
as a graded vector space.
\qed

\sssn
\Lemma
Let $M$ be a graded $\n$-module with a filtration
$$
F^0M\subset\ldots\subset F^mM\subset\ldots\subset M,\
\dim F^mM<\infty\text{ for every }m\in\Z_{>0}.
$$
Then
$
\tor_{\ne\dim\n}^{\n}(U(\n)^*,M)=0,\
\tor_{\dim\n}^{\n}(U(\n)^*,M)=M.
$

\dok
Recall that the functor $\underset{\map}{\lim}$ is exact. The standard
complex for the computation of $\n$-homology
$$
U(\n)^*\ten\Lambda(\n)\ten M
=\underset{\map m}{\lim }
U(\n)^*\ten\Lambda(\n)\ten F^mM.
$$
Thus we have
$$
H\bul(U(\n)^*\ten\Lambda(\n)\ten M)
=\underset{\map m}{\lim }H\bul
(U(\n)^*\ten\Lambda(\n)\ten F^mM)=\delta_{\bullet,-\dim\n}\underset{\map m}
{\lim}F^mM=\delta_{\bullet,-\dim\n}M. \qed
$$
Now recall that $L(\lambda)$ is an integrable $\g$-module, in particular it
is a union of finite dimensional $\n_w$-modules. Using the previous Lemma we
see that the statement of the Proposition is proved on the level of
$\h^*\times\Z$-graded vector spaces. But an integrable $\g$-module is
completely determined by its character.  \qed

\sssn
\Cor
There exists a complex $B_w\bul(\lambda)$ in $\M_k$ of the form
\begin{multline*}
\ldots
\map
\underset{{v\in W,}\atop{\lth(v)=m}}{\bigoplus}
M_w(wv\cdot\lambda)\langle-\hgt(\lambda-wv\cdot\lambda)\rangle
\map\ldots\\ \map
\underset{{v\in W,}\atop{\lth(v)=1}}{\bigoplus}
M_w(wv\cdot\lambda)\langle-\hgt(\lambda-wv\cdot\lambda)\rangle
\map
M_w(w\cdot\lambda)\langle -\hgt(\lambda-w\cdot\lambda)\rangle\map 0
\end{multline*}
such that
$H^{\ne-\lth(w)}(B_w\bul(\lambda))=0,\
H^{-\lth(w)}(B_w\bul(\lambda))=L(\lambda).
$
Here as usual $\langle\cdot\rangle$ denotes the shift of the $\Z$-grading.
\qed

\sssn
Recall that in~\ref{tlength} we have defined the twisted length function on
the affine Weyl group with the twist $w\in W$ by
$$
\lth^w(u)=\lth(w^{-1}u)-\lth(w^{-1}).
$$
\sssn
\Cor
The complex $B_w(\lambda)[-\lth(w)]$ can be rewritten as follows:
\begin{multline*}
\ldots
\map
\underset{{v\in W,}\atop{\lth^w(v)=m}}{\bigoplus}
M_w(v\cdot\lambda)\langle-\hgt(\lambda-v\cdot\lambda)\rangle
\map\ldots \map
\underset{{v\in W,}\atop{\lth^w(v)=0}}{\bigoplus}
M_w(v\cdot\lambda)\langle-\hgt(\lambda-v\cdot\lambda)\rangle
\map\ldots\\
\map
\underset{{v\in W,}\atop{\lth^w(v)=-\lth(w)+1}}{\bigoplus}
M_w(v\cdot\lambda)\langle-\hgt(\lambda-v\cdot\lambda)\rangle
\map
M_w(w\cdot\lambda)\langle -\hgt(\lambda-w\cdot\lambda)\rangle\map 0.\qed
\end{multline*}
\sssn
\Def
We call  the complex $B_w\bul(\lambda)$ the {\em twisted BGG
resolution} of the module $L(\lambda)$ with the twist $w$.

\section{Limit procedure and semi-infinite BGG resolution}

\subsection{Recollections of semi-infinite Lie algebra cohomology.}
Here we present a list of the definitions and facts about semi-infinite
cohomology of Lie algebras to  be used later.

\subsubsection{Critical cocycle.}
Let $V=\underset {n\in\Z}\bigoplus V_n$ be a graded vector
space such that $\dim V_n<\infty$ for $n>0$. Denote the
space $\underset{n\le 0}{\bigoplus}V_n$ (resp. the space
$\underset{n>0}{\bigoplus}V_n$) by $V^-$ (resp. by $V^+$).
Consider the Lie algebra $\gl(V)$ of linear transformations
of $V$ satisfying the condition:
$$
 a\in\gl(V),\
 a=(a_{ij}), \Longrightarrow \exists m\in\N:\ a_{ij}=0 \text{
 for }|i-j|>m.
$$
$\gl(V)$ has a well known central
extension given by the 2-cocycle $\omega_0$:
$$
 \omega_0(a_1,a_2)=\tr(\pi_{V^+}\circ[a_1,a_2]-[\pi_{V^+}\circ
 a_1,\pi_{V^+}\circ a_2]).
$$
Here $\pi_{V^+}$ denotes the projection
$V\map V_+$ with the kernel $V_-$. Note that the trace is
well defined on maps $V_+\map V_+$ nontrivial on a finite number of the
grading components.

\sssn
Let  $\a$ be a graded Lie algebra such that all $\a_m$ are finite dimensional.
Then $\a=\ap\oplus\am$ as a vector space, where
$$
\ap:=\underset{m>0}{\bigoplus}\a_m,\
\am:=\underset{m\le0}{\bigoplus}\a_m.
$$
For any graded $\a$-module $V$ we have a morphism of Lie algebras
$\a\map\gl(V)$. The inverse image of the 2-cocycle $\omega_0$ is denoted by
$\omega_0^V$. In particular the adjoint representation of $\a$ provides the
2-cocycle $\omega_0^{\a}$ called the {\em critical cocycle} of $\a$.

\subsubsection{Feigin standard complex.}
Choose a homogenious base $\{e_i|i\in\Z\}$ in $\a$ and let
$\{c_{i,j}^k\}_{i,j,k\in\Z}$ denote the structure constants of $\a$ in this
base.  For a graded $\a$-module $M=\underset{m\le m_0}{\bigoplus}M_m$
consider a graded vector space
$$
\stand(M):=\Lambda\bul(\a^{+*})\ten\Lambda\bul(\am)\ten M,\
\deg \overe_n^*=1,\ \deg \overe_m=-1,\ \overe_m\in\am,\ \overe_n^*\in\ap.
$$
Then the element
$
D:=\underset{i\in\Z}{\sum} \overe_i^*\ten e_i+\underset{i<j,k}{\sum}
c_{i,j}^k:\overe_i^*\overe_j^*\overe_k:
$
can be considered as an endomorphism of $\stand(M)$.
Here as usual $:\ :$ denotes the normal ordering.  It is known that the
cohomology class of the  critical 2-cocycle $\omega_0^{\a}$ of the Lie
algebra $\a$ is the only obstruction for $D^2=0$. We  suppose below that the
cohomology class of $\omega_0^{\a}$ equals to zero and $D^2=0$. Then
$\stand(M)$ becomes a complex of vector spaces with the differential $D$
called  the {\em Feigin standard complex}.

\sssn              \label{cond}
\Rem
Assume that $\a$ has a decreasing filtration by subalgebras
$$
\a=F^0\a\supset\ldots\supset F^m\a\supset\ldots,\ [F^m\a,F^n\a]\subset
F^{m+n}\a
$$
and all the inclusions $F^{m+1}\a\subset F^m\a$ are strict.  In
this case the 2-cocycle $\omega_0$ itself equals to zero.

\sssn
\Def
The cohomologies  of the complex $\stand(M)$ are called {\em semi-infinite
cohomology spaces} of the Lie algebra $\a$ with coefficients in the module $M$:
$$
H^{\si+\bullet}(M):=H\bul(\stand(M)).
$$
\sssn \label{limmmap}
Let $\a=\underset{m\in\Z}{\bigoplus}\a_m$ be a graded Lie algebra containing
two graded Lie subalgebras $\b_1$ and $\b_2$ such that $\a_{> 0}=\b_{1> 0}$
and $\a_{\le0}=\b_{2\le0}$.  We suppose that the condition~\ref{cond} is
satisfied and the obstruction 2-cocycle is zero for all three Lie algebras.

\sssn
\Lemma                  \label{maps}
For a graded $\a$-module $M=\underset{m<m_0}{\bigoplus}M_m$ there exist
natural morphisms
$$
H^{\si+\bullet}(\a,M)\map
H^{\si+\bullet}(\b_2,M)\text{ and }
H^{\si+\bullet}(\b_1,M)\map
H^{\si+\bullet}(\a,M).
$$

\dok
There exist natural morphisms of the graded vector spaces
$$
\stand(\a,M)\map
\stand(\b_2,M)\text{ and }
\stand(\b_1,M)\map
\stand(\a,M).
$$
One has to check that the maps commute with the differentials and provide
morphisms of the standard complexes for the computation of Lie algebra
semi-infinite cohomology.  \qed

\sssn
\Rem             \label{remm}
The statement of the previous Lemma follows easily (and perhaps more
naturally) from the general homological algebra (see [V] 3.9, [Ar1] Appendix
B). Namely semi-infinite cohomology of a Lie algebra $\a$
satisfying~\ref{cond} can be viewed as an exotic derived functor of a functor
$M\mapsto\hom_{\a}(\underc,S_{\a_{\le0}}^{\a}\ten_{U(\a)} M)$. Now in the
Lemma we have
$$
S_{\a_{\le0}}^{\a}
=\Ind_{\b_2}^{\a}
S_{\b_{2\le0}}^{\b_2},\text{ and }
S_{\a_{\le0}}^{\a}
=\Coind_{\b_1}^{\a}
S_{\b_{1\le0}}^{\b_1}
$$
Thus there exist natural morphisms
\begin{gather*}
\hom_{\a}(\underc,S_{\a_{\le0}}^{\a}\ten_{U(\a)} M)\map
\hom_{\b_2}(\underc,S_{\a_{\le0}}^{\a}\ten_{U(\a)} M)\til{\map}
\hom_{\b_2}(\underc,S_{\b_{2\le0}}^{\b_2}\ten_{U(\b_2)} M),\\
\hom_{\b_1}(\underc,S_{\b_{1\le0}}^{\b_1}\ten_{U(\b_1)} M)\til{\map}
\hom_{\a}(\underc,S_{\a_{\le0}}^{\a}\ten_{U(\b_1)} M)\map
\hom_{\a}(\underc,S_{\a_{\le0}}^{\a}\ten_{U(\a)} M)
\end{gather*}
that provide the required morphisms of the derived functors.

\subsubsection{Semi-infinite standard resolutions.}
Let $\a=\underset{m\in\Z}{\bigoplus}\a_m$ be a graded Lie algebra
with a graded subalgebra $\b$ such that $\b_{>0}$ is
of finite codimension in $\g_{>0}$ and the 2-cocycle $\omega_0^{\b}=0$.  Here
we construct a semi-infinite analogue of the standard resolution~\ref{stres}.
Recall the following statement crucial for Lie algebra semi-infinite
cohomology.

\sssn \label{end}
\Prop
(see [V] 3.2.1, [Ar2] Corollary 4.4.2)
Let $\til{\a}=\a\oplus\CC K$ be the central extension of $\a$ with the help
of the cocycle $\omega_0^{\a}$. Let $\til{U(\a)}:=U(\til{\a})/\{K-1\}$.  Then
there exists an inclusion of associative algebras
$\til{U(\a)}\hookrightarrow\End_{\a}(S_{\a_{>0}}^{\a})$.  \qed

We suppose below that the inclusion $\b\subset\a$ can be extended to
$\b\subset\til{U(\a)}$.

For a graded vector space $V=\underset{m\in\Z}{\bigoplus}V_m$ consider
the semi-infinite exterior powers
$$
\Lambda^{\si+\bullet}(V):=\Lambda(V_{>0}^*)\ten\Lambda(V_{\le0}).
$$
\Lemma
(see e.~g.~[V] Proposition 2.3)
Let $V$ be a graded representation of a graded Lie algebra $\a$. Then
$\Lambda^{\si+\bullet}(V)$ is a module over the central extension of
$\a$ with the help of $\omega_0^V$.
\qed

In particular the semi-infinite exterior
powers of the $\b$-module $\a/\b$ are $\b$-modules as well.

\sssn  \label{filt}
Consider the Feigin standard complex
$\stand(S_{\a_{>0}}^{\a})=S_{\a_{>0}}^{\a}\ten
\Lambda^{\si+\bullet}(\a)\ten\underc$ for the computation of the
semi-infinite Lie algebra cohomology of $\a$ with coefficients in the
semiregular $\a$-module (using the second $\a$-module structure on
$S_{\a_{>0}}^{\a}$ from Proposition~\ref{end}).  Clearly we have
$$
H^{\si+\bullet}(S_{\a_{>0}}^{\a})=\delta_{\bullet,0}\underc
$$
and $\stand(S_{\a_{>0}}^{\a})\cong
S_{\a_{>0}}^{\a}
\ten\Lambda^{\si+\bullet}(\b)\ten\Lambda^{\si+\bullet}(\a/\b)$ as a vector
space.  We introduce a filtration of the complex as follows:
$$
F^m(\stand(S_{\a_{>0}}^{\a})):=\underset{{p+q=n}\atop{q\ge m}}{\bigoplus}
S_{\a_{>0}}^{\a}\ten\Lambda^{\si+p}(\b)\ten\Lambda^{\si+q}(\a/\b).
$$
Clearly the differential in the complex preserves the filtration. Consider
the corresponding spectral sequence. We have
$$
E_1^{p,q}=H^{\si+p}(\b,S_{\a_{>0}}^{\a}\ten\Lambda^{\si+q}(\a/\b))=
\delta_{0,p}H^{\si+0}(\b,S_{\a_{>0}}^{\a}\ten\Lambda^{\si+q}(\a/\b)).
$$
Thus the spectral sequence converges in $E_2$. On the other hand
$E_{\infty}^{p,q}=\delta_{p,0}\delta_{q,0}\underc$ and we obtain
the following statement.

\sssn
\Lemma
There exists a complex of $\a$-modules
$$
\St_{\b}:=
H^{\si+0}(\b,S_{\a_{>0}}^{\a}\ten
\Lambda^{\si+\bullet}(\a/\b))
$$
such that $H^{\ne0}(\St_{\b})=0,\
H^{0}(\St_{\b})=\underc$.
\qed

We call the obtained complex the {\em semi-infinite standard resolution} of
the trivial $\a$-module relative to the subalgebra $\b$.

\sssn \label{stress}
For a $\g$-module $M$ consider the complex of $\g$-modules $\St_{\b}\ten M$
(where the $\g$-module structure uses the comultiplication on $U(\g)$).
Evidently it is quasiisomorphic to $M$.
We call the obtained complex the standard resolution of $M$ and denote it by
$\St_{\b}(M)$.

Let $\b_1$, $\b_2$ and $\b$ be subalgebras in $\a$ such that
$$
\b\subset\b_1,\ \b\subset \b_2,\ \b_{>0}=\b_{1>0},\
\b_{\le0}=\b_{2\le0}\text{ and }
\omega_0^{\b_1}=\omega_0^{\b_2}=\omega_0^{\b}=0.
$$
\sssn
\Lemma \label{morph}
There exist natural quasiisomorphisms
$$
\St_{\b_2}(M)\map\St_{\b}(M)\map\St_{\b_1}(M).
$$
\dok
Recall that semi-infinite standard resolutions were constructed using levels
$E_1$ of certain spectral sequences (see~\ref{filt}).
Since $\b\subset\b_1$ and $\b_{>0}=\b_{1>0}$ we have
\begin{gather*}
F^m_{\b}(\stand(S_{\a_{>0}}^{\a}))
=\underset{{{p+q=n}\atop{q\ge m}}}{\bigoplus}
S_{\a_{>0}}^{\a}\ten\Lambda^{\si+p}(\b)\ten\Lambda^{\si+q}(\a/\b)\\
=\underset{{p+q=n}\atop{q\ge m}}{\bigoplus}
S_{\a_{>0}}^{\a}\ten
\underset{m+n=p}{\bigoplus}
\Lambda^{-m}(\b_{\le0})\ten
\Lambda^{n}(\b_{>0}^*)\ten
\underset{s+t=q}{\bigoplus}
\Lambda^{-s}((\a/\b)_{\le0})\ten
\Lambda^{t}((\a/\b)_{>0}^*) \\ \subset
\underset{{p+q=n}\atop{q\ge m}}{\bigoplus}
S_{\a_{>0}}^{\a}\ten
\underset{m+n=p}{\bigoplus}
\Lambda^{-m}(\b_{1\le0})\ten
\Lambda^{n}(\b_{1>0}^*)\ten
\underset{s+t=q}{\bigoplus}
\Lambda^{-s}((\a/\b_1)_{\le0})\ten
\Lambda^{t}((\a/\b_1)_{>0}^*) \\
=
F^m_{\b_1}(\stand(S_{\a_{>0}}^{\a})).
\end{gather*}
Since $\b\subset\b_2$ and $\b_{\le0}=\b_{2\le0}$ we have
\begin{gather*}
F^m_{\b_2}(\stand(S_{\a_{>0}}^{\a}))
=\underset{{p+q=n}\atop{q\ge m}}{\bigoplus}
S_{\a_{>0}}^{\a}\ten\Lambda^{\si+p}(\b_2)\ten\Lambda^{\si+q}(\a/\b_2) \\
=\underset{{p+q=n}\atop{q\ge m}}{\bigoplus}
S_{\a_{>0}}^{\a}\ten
\underset{m+n=p}{\bigoplus}
\Lambda^{-m}(\b_{2\le0})\ten
\Lambda^{n}(\b_{2>0}^*)\ten
\underset{s+t=q}{\bigoplus}
\Lambda^{-s}((\a/\b_2)_{\le0})\ten
\Lambda^{t}((\a/\b_2)_{>0}^*) \\ \subset
\underset{{p+q=n}\atop{q\ge m}}{\bigoplus}
S_{\a_{>0}}^{\a}\ten
\underset{m+n=p}{\bigoplus}
\Lambda^{-m}(\b_{\le0})\ten
\Lambda^{n}(\b_{>0}^*)\ten
\underset{s+t=q}{\bigoplus}
\Lambda^{-s}((\a/\b)_{\le0})\ten
\Lambda^{t}((\a/\b)_{>0}^*)    \\
=
F^m_{\b}(\stand(S_{\a_{>0}}^{\a})).
\end{gather*}
Thus we obtain the required morphisms of the spectral sequences that
correspond to the morphisms of the filtered complexes.  \qed

\subsection{Twisted Verma modules and Wakimoto modules.}
Consider the Lie algebra $\dot w(\b^+)$ with the $\Z$-grading induced by the
one on $\g$ (see~\ref{gr}). Here we give a description of twisted Verma
modules in terms of semi-infinite cohomology of $\dot w(\b^+)$.

Consider the semiregular $U_k(\g)$-module with respect to $\n^+$.
The following statement is proved in  [Ar2], 4.5.9. It is in fact a particular
case of Proposition~\ref{end}.

\sssn
\Prop
$\End_{U_k(\g)}(S_{\n^+})\supset U_{-2h^\vee-k}(\g)$, where $h^\vee$
denotes the dual Coxeter number of $\overg$.
\qed

Consider $S_{\n^+}^{\g}$ as a  $\dot w(\b^+)$-module using this second
$\g$-module structure.

\sssn
\Lemma
$$
H^{\frac \infty 2+m}(\dot w(\b^+),S_{\n^+}^{\g}\ten\CC(\lambda))=0
\text{ when }m\ne 0,\
H^{\frac \infty 2+0}(\dot w(\b^+),S_{\n^+}^{\g}\ten\CC(\lambda))
=S_{\dot w(\n_w)}^{\dot w(\n^-)}\ten\CC(\lambda)
$$
as a $\h^*\times\Z$-graded vector space.
\qed

Thus we have a $U_k(\g)$-module $\til{M}_w(\lambda):=
H^{\frac \infty 2+0}(\dot w(\b^+),S_{\n^+}^{\g}\ten\CC(\lambda))$
that has the size of the twisted Verma module.

\sssn
\Rem \label{cofr}
By definition of the $\g$-module structure on $S_{\n^+}^{\g}$ the module
$\til{M}_w(\lambda)$ is cofree over $U(\dot w(\n_w))$.

\sssn
Consider also the Lie algebra $\b^{\si}:=\overh\ten\CC[t]
\oplus\overg^-\ten \CC[t,t^{-1}]$ and a $\h^*\times\Z$-graded
$U_k(\g)$-module $W(\lambda):=
H^{\frac \infty 2+0}(\b^{\si},S_{\n^+}^{\g}\ten\CC(\lambda))$.

\Def
The $\g$-module $W(\lambda)$ is called the {\em Wakimoto module} at the level
$k$ of the weight $\lambda$.

We show that the modules $\til{M}_w(\lambda)$ are indeed isomorphic to
twisted Verma modules with the twist $w$ and construct morphisms between
them so that $W(\lambda)$ becomes their inductive limit as $\lth(w)$ tends
to infinity.

\subsection{Equivalence of definitions of twisted Verma modules.}
Recall that a twisted Verma module with the highest weight $\lambda$ is
defined as $\Phi_w(M(w^{-1}\cdot
\lambda))\langle\hgt(w^{-1}\cdot\lambda-\lambda)\rangle$.

\sssn
\Lemma  \label{e}
$\til{M}_w(\lambda) \cong M_w(\lambda)$.

\dok
Up to a $\Z$-grading shift we have
\begin{gather*}
\Phi_w(M(w^{-1}\cdot\lambda))=T_w\circ S_w
(H^0(\b^+,S_{\n^+}^{\g}\ten\CC(w^{-1}\cdot\lambda)))\\ =H^0(\dot w(\b^+),T_w(
(\Coind_{\dot w^{-1}(\n^-)\cap\n^-\oplus\h}^{\g} \Ind_{\CC}^{\dot
w^{-1}(\n^-)\cap\n^-\oplus\h}\underc) \ten\CC(w^{-1}\cdot\lambda)))\\
=H^0(\dot w(\b^+),S_{\dot w(\b^+)\cap\n^-}^{\g}\ten_{U_k(\g)}(S_{\n^+}^{\g}
\ten\CC(\lambda))) =H^0(\dot w(\b^+),S_{\dot w(\b^+)\cap\n^-}^{\dot w(\b^+)}
\ten_{U(\dot w(\b^+))}(S_{\n^+}^{\g} \ten\CC(\lambda))).
\end{gather*}
Now recall the abstract homological definition of Lie algebra semi-infinite
cohomology from Remark~\ref{remm}. The $\dot w(\b^+)$-module
$S_{\n^+}^{\g}\ten\CC(\lambda)$ is {\em semijective} (see [V] 3.1, [Ar1]
Appendix B for the definition of semijectiveness).  Thus we need no
resolution to find the 2-sided derived functor:
$$
H^0(\dot w(\b^+),S_{\dot w(\b^+)\cap\n^-}^{\dot w(\b^+)}
\ten_{U(\dot w(\b^+))}(S_{\n^+}^{\g}
\ten\CC(\lambda)))
=H^{\si+0}(\dot w(b^+),S_{\n^+}^{\g}\ten\CC(\lambda)).\qed
$$
\sssn
\Prop
$$
\St_{w(\b^+)}(L(\lambda))\cong
\Phi
_{w}(D\bul(L(\lambda))
[-\lth(w)]
$$
where $[\cdot]$ denotes the shift of the homological grading.

\dok
The statement follows immediately from Lemmas~\ref{1}
and~\ref{2}.
\qed

\sssn \label{1}
\Lemma
$\St_{\b^+}\cong D\bul$.
\qed

\sssn
\Lemma \label{2}
$\Phi_w(\St_{\b^+})\cong \St_{\dot w(\b^+)}[\lth(w)]$.

\dok
The proof is parallel to the one of Lemma~\ref{e}.
\qed

The next statement shows that the twisted BGG resolutions are obtained from
the semi-infinite standard resolutions relative to $\dot w(\b^+)$ in the same
way as the ordinary BGG resolution was obtained from the ordinary standard
resolution relative to $\b^+$ (see~\ref{stres}).

\sssn
\Lemma
For every $m\in\Z$ the $\g$-module
$\St_{\dot w(\b^+)}(L(\lambda))$ is finitely filtered by twisted Verma
modules with the twist $w$.  \qed

\sssn
\Cor   \label{cut}
$
B_w\bul(\lambda)\cong
\left(\St_{\dot w(\b^+)}(L(\lambda))
\right)_{(\lambda+\rho,\lambda+\rho)}[\lth(w)]
$.
\qed

\subsection{From twisted Verma modules to Wakimoto modules.}
Fix $w_1,w_2\in W$ such that $\lth(w_2w_1)=\lth(w_1)+\lth(w_2)$.
Consider the following Lie subalgebras in $\g$:
$$
\b_{w_2w_1,w_1}:=\dot w_1(\b^+)\cap\dot w_2\dot w_1(\b^+),\
\b_{w_2w_1,w_1}^{w_1}:=
\underset{\alpha\in R_{w_2w_1,w_1}^{w_1}}{\bigoplus}\g_{\alpha},\
\b_{w_2w_1,w_1}^{w_2w_1}:=
\underset{\alpha\in R_{w_2w_1,w_1}^{w_2w_1}}{\bigoplus}\g_{\alpha},
$$
where
$R_{w_2w_1,w_1}^{w_1}:=\{\alpha\in w_1(R^+)|w_2(\alpha)\not\in w_1(R^+)\},\
R_{w_2w_1,w_1}^{w_2w_1}:=\{\alpha\in w_2w_1(R^+)|w_1^{-1}(\alpha)\not\in
w_1(R^+)\} $.

Then
$
\b_{w_2w_1,w_1}\oplus
\b_{w_2w_1,w_1}^{w_1}=\dot w_1(\b),
\b_{w_2w_1,w_1}\oplus
\b_{w_2w_1,w_1}^{w_2w_1}=\dot w_2\dot w_1(\b)$.
By the choice of $w_1$ and $w_2$
we have $\b_{w_2w_1,w_1}^{w_1}\subset\n^+$ and
$\b_{w_2w_1,w_1}^{w_2w_1}\subset\n^-$.
In particular using~\ref{maps} one can define morphisms of $\g$-modules
$$
H^{\si+\bullet}(\dot w_1(\b^+),S_{\n^+}^{\g}\ten\CC(\lambda))\map
H^{\si+\bullet}(\b_{w_2w_1,w_1},S_{\n^+}^{\g}\ten\CC(\lambda))\map
H^{\si+\bullet}(\dot w_2\dot w_1(\b^+),S_{\n^+}^{\g}\ten\CC(\lambda)),
$$
i.~e. $M_{w_1}(\lambda)\map M_{w_2w_1}(\lambda)$. Denote these morphisms
by $\varphi_{w_2w_1,w_1}^{\lambda}$.

For $w\in W$ such that $\dot w(\b^+)_{\le0}\cap\b^{\si}=
\dot w(\b^+)_{\le0}$ and $\dot w(\b^+)\cap\b_{>0}^{\si}=\b^{\si}$
consider  maps

$$
H^{\si+\bullet}(\dot w(\b^+),S_{\n^+}^{\g}\ten\CC(\lambda))\map
H^{\si+\bullet}(\dot w(\b^+)\cap\b^{\frac\infty
2},S_{\n^+}^{\g}\ten\CC(\lambda))\map
H^{\si+\bullet}(\b^{\si},S_{\n^+}^{\g}\ten\CC(\lambda)),
$$
i.~e. $M_w(\lambda)\map W(\lambda)$.
Denote these morphisms
by $\varphi_{\si,w}^{\lambda}$.

\sssn  \label{w}
\Lemma

\qquad(i) For $w_1,w_2,w_3\in W$ such that $\lth(w_3w_2w_1)=\lth(w_1)+
\lth(w_2)+\lth(w_3)$ we have
$$
\varphi_{w_3w_2w_1,w_2w_1}^{\lambda}\circ
\varphi_{w_2w_1,w_1}^{\lambda}=
\varphi_{w_3w_2w_1,w_1}^{\lambda}.
$$
\qquad(ii)
For $w_1,w_2\in W$ such that $\lth(w_2w_1)=\lth(w_1)+
\lth(w_2)$ we have
$
\varphi_{\si,w_2w_1}^{\lambda}\circ
\varphi_{w_2w_1,w_1}^{\lambda}=
\varphi_{\si,w_1}^{\lambda}.
$
\qed

Fix $w_1,w_2\ldots\in W$ such that for each $m$ we have
$\lth(w_m\ldots w_1)=\lth(w_1)+\ldots+\lth(w_m)$.
Thus for each $\lambda$ we obtain an inductive system consisting of twisted
Verma modules with twists $w_1, w_2w_1,\dots\in W$ and the lengths of
$w_m\ldots w_1$ tending to infinity. Suppose the chosen sequence of elements
of the Weyl group belongs to the image of $-Q^{\prime\prime+}$ in $T$. Then
we obtain the following statement.

\sssn
\Lemma
$W(\lambda)=\underset{\map}{\lim}M_{w_m\ldots w_1}(\lambda)$.

\dok
The maps $\varphi_{w_m\ldots w_1,w_{m-k}\ldots w_1}^{\lambda}$
are well defined as maps
$$
\stand(\dot w_{m-k}\ldots\dot w_1(\b^+),S_{\n^+}^{\g}\ten\CC(\lambda))\map
\stand(\dot w_m\ldots\dot w_1(\b^+),S_{\n^+}^{\g}\ten\CC(\lambda)).
$$
On the other hand we have
$$
\stand(\b^{\si},S_{\n^+}^{\g}\ten\CC(\lambda))=
\underset{\map}{\lim}
\stand(\dot w_m\ldots\dot w_1(\b^+),S_{\n^+}^{\g}\ten\CC(\lambda)).
$$
But the direct image functor is exact.
\qed

We are almost done. It remains to prove that roughly speaking
the marphisms $\varphi_{w_2w_1,w_1}^{\lambda}$
indeed form an inductive system of complexes
(twisted BGG resolutions).
The combinalorial proof of this statement requires the following fact:
$$
\dim\hom_{\g}(M(\lambda),M_w(\mu))\le1.
$$
Unfortunately I do not know any proof of the fact, and our considerations
become more cumbersome.

\subsection{The limit procedure.}
Fix the sequence $w_1,\ldots,w_m,\ldots\in W$
satisfying the conditions from the previous paragraph.
Consider the standard resolutions
$\St_{\dot w_m\ldots \dot w_1(\b^+)}(L(\lambda))$ (see~\ref{stress}).
Then for every $m$ and $k$ the algebras
$\dot w_{m+k}\ldots\dot w_1(\b^+)$,
$\dot w_m\ldots\dot w_1(\b^+)$
and
$\dot w_{m+k}\ldots\dot w_1(\b^+)\cap
\dot w_m\ldots\dot w_1(\b^+)$ satisfy the conditions of  Lemma~\ref{morph},
and we obtain canonical morphisms of complexes
$$
\psi_{
w_{m+k}\ldots w_1,
w_m\ldots w_1}^{\lambda}:\
\St_{
\dot w_m\ldots\dot w_1(\b^+)}(M)\map
\St_{
\dot w_{m+k}\ldots\dot w_1(\b^+)\cap
\dot w_m\ldots\dot w_1(\b^+)}(M)\map\St_{
\dot w_{m+k}\ldots\dot w_1(\b^+)}(M).
$$
Thus we obtain an inductive system complexes, and all morphisms of complexes
in it are quasiisomorphisms by~\ref{morph}.

\sssn
\Theorem
There exists a complex in $\M_k$ quasiisomorphic to $L(\lambda)$ of the form
\begin{multline*}
\ldots
\map
\underset{{w\in W,}\atop{\lth^{\si}(w)=m}}{\bigoplus}
W(w\cdot\lambda)\langle-\hgt(\lambda-w\cdot\lambda)\rangle
\map\ldots\map
\underset{{w\in W,}\atop{\lth^{\si}(w)=0}}{\bigoplus}
W(w\cdot\lambda)\langle-\hgt(\lambda-w\cdot\lambda)\rangle
\map            \\
\underset{{w\in W,}\atop{\lth^{\si}(w)=-1}}{\bigoplus}
W(w\cdot\lambda)\langle-\hgt(\lambda-w\cdot\lambda)\rangle\map
\ldots
\underset{{w\in W,}\atop{\lth^{\si}(w)=-m}}{\bigoplus}
W(w\cdot\lambda)\langle-\hgt(\lambda-w\cdot\lambda)\rangle\map
\ldots
\end{multline*}
\dok
We know already from Corollary~\ref{cut} that when restricted to
the $\Gamma$-homogenious direct summands with  eigenvalue of $\Gamma$ equal to
$(\lambda+\rho,\lambda+\rho)$ the morphisms
$$
\psi_{
w_{m+k}\ldots w_1,
w_m\ldots w_1}^{\lambda}:\
\St_{
\dot w_m\ldots\dot w_1(\b^+)}(L(\lambda))\map
\St_{\dot w_{m+k}\ldots\dot w_1(\b^+)}(L(\lambda))
$$
form an inductive system of morphisms between twisted BGG resolutions
$$
\psi_{
w_{m+k}\ldots w_1,
w_m\ldots w_1}^{\lambda}:\
\B\bul_{w_m\ldots w_1}(\lambda)\map
\B_{w_{m+k}\ldots w_1}\bul(\lambda)[-\lth(w_{m+1})-\ldots-\lth(w_{m+k})].
$$
Now we choose a sequence of elements of the affine Weyl group satisfying
the conditions from~\ref{w} and belonging to the image of $-Q^{\prime\prime+}$
in $T\subset W$.  Thus by Lemma~\ref{maincomb} for any $w\in W$ there exists
$m_0\in\N$ such that for any $m\ge m_0$
$$
M_{w_m\ldots w_1}(w\cdot\lambda)\langle-\hgt(\lambda-w\cdot\lambda)\rangle
\hookrightarrow B^{\lth^{\si}(w)}_{w_m\ldots w_1}(\lambda)
[-\lth(w_m)-\ldots-\lth(w_1)].
$$
It remains to check that for every $w\in W$ and $m$  sufficiently large
we have
$$
\psi_{w_{m+k}\ldots w_1,
w_m\ldots w_1}^{\lambda}
|_{M_{w_m\ldots w_1}(w\cdot\lambda)}=
\varphi_{w_{m+k}\ldots w_1,
w_m\ldots w_1}^{w\cdot\lambda}.
$$
\qed

\section*{References}
$\text{[Ar1]}$ S.Arkhipov. {\it Semiinfinite cohomology of quantum groups.}
Preprint q-alg/9601026 (1996), 1-24.\\ $\text{[Ar2]}$ S.Arkhipov. {\it
Semiinfinite cohomology of associative algebras and bar duality.} Preprint
q-alg/9602013 (1996), 1-21.\\ $\text{[BGG]}$ I.~I.~Bernstein, I.~M.~Gelfand,
S.~I.~Gelfand.  {\it Differential operators on the principal affine space and
investigation of} $\g${\em -modules}, in Proceedings of Summer School on Lie
groups of Bolyai J\=anos Math. Soc., Helstead, New York, 1975.\\
$\text{[FF]}$ B. Feigin, E. Frenkel. {\em Affine Kac-Moody algebras and
semi-infinite flag manifolds.} Comm. Math. Phys. \hbox{\bf 128} (1990) ,
161-189.\\ $\text{[K]}$ V. Kac. {\it Infinite dimensional Lie algebras.}
Birkh\"auser, Boston (1984).\\ $\text{[L1]}$ G. Lusztig. {\it Monodromic
systems on affine flag manifolds.} Proc.R. Soc. Lond. A \hbox{\bf 445} (1994)
231-246.\\ $\text{[L2]}$ G. Lusztig. {\it Hecke algebras and Jantzen's
generic decomposition patterns.} Adv. in Math. Vol.\hbox{\bf 37},
No.\hbox{\bf 2} (1980) 121-164.\\ $\text{[RW]}$ A.~Rocha-Charidi,
N.~R.~Wallach. {\em Projective modules over graded Lie algebras.}
Math.~Z.~\hbox{\bf 180}, (1982), 151-177.\\ $\text{[V]}$ A.Voronov. {\it
Semi-infinite homological algebra.} Invent. Math. {\bf 113}, (1993),
103--146.
\end{document}